
\def\LaTeX{%
  \let\Begin\begin
  \let\End\end
  \let\salta\relax
  \let\finqui\relax
  \let\futuro\relax}

\def\UK{\def\our{our}\let\sz s\def\analogue{analogue}}
\def\USA{\def\our{or}\let\sz z\def\analogue{analog}}



\LaTeX

\USA


\salta

\documentclass[twoside,a4paper,12pt]{article}
\setlength{\textheight}{24cm}
\setlength{\textwidth}{16cm}
\setlength{\oddsidemargin}{2mm}
\setlength{\evensidemargin}{2mm}
\setlength{\topmargin}{-15mm}
\parskip2mm


\usepackage[usenames,dvipsnames]{color}
\usepackage{amsmath}
\usepackage{amsthm}
\usepackage{amssymb}
\usepackage[mathcal]{euscript}
%
%
\usepackage{cite}
%
%
%


\definecolor{viola}{rgb}{0.3,0,0.7}
\definecolor{ciclamino}{rgb}{1,0,1}
\definecolor{rosso}{rgb}{0.8,0,0}

\def\gianni #1{{\color{rosso}#1}}
\def\pier #1{{\color{rosso}#1}}

\def\betti #1{{\color{ciclamino}#1}}
\def\rev #1{{\color{rosso}#1}}
\def\revdue #1{{\color{blue}#1}}


\def\betti #1{#1}
\def\gianniold #1{#1}
\def\gianni #1{#1}
\def\pier #1{#1}
\def\rev #1{{#1}}
\def\revdue #1{{#1}}



\bibliographystyle{plain}


%

\finqui

\def\Beq{\Begin{equation}}
\def\Eeq{\End{equation}}
\def\Bsist{\Begin{eqnarray}}
\def\Esist{\End{eqnarray}}

\def\Bthm{\Begin{theorem}}
\def\Ethm{\End{theorem}}
\def\Blem{\Begin{lemma}}
\def\Elem{\End{lemma}}

\def\Brem{\Begin{remark}\rm}
\def\Erem{\End{remark}}

\def\Bdim{\Begin{proof}}
\def\Edim{\End{proof}}
\def\Bcenter{\Begin{center}}
\def\Ecenter{\End{center}}
\let\non\nonumber




\def\step #1 \par{\medskip\noindent{\bf #1.}\quad}


\def\Lip{Lip\-schitz}
\def\Holder{H\"older}

\def\aand{\quad\hbox{and}\quad}

\def\lhs{left-hand side}
\def\rhs{right-hand side}
\def\sfw{straightforward}


\def\organiz{organi\sz}
\def\recogniz{recogni\sz}
\def\regulariz{regulari\sz}

\def\bhv{behavi\our}


\def\multibold #1{\def\arg{#1}%
  \ifx\arg\pto \let\next\relax
  \else
  \def\next{\expandafter
    \def\csname #1#1#1\endcsname{{\bf #1}}%
    \multibold}%
  \fi \next}

\def\pto{.}

\def\multical #1{\def\arg{#1}%
  \ifx\arg\pto \let\next\relax
  \else
  \def\next{\expandafter
    \def\csname cal#1\endcsname{{\cal #1}}%
    \multical}%
  \fi \next}


\def\multimathop #1 {\def\arg{#1}%
  \ifx\arg\pto \let\next\relax
  \else
  \def\next{\expandafter
    \def\csname #1\endcsname{\mathop{\rm #1}\nolimits}%
    \multimathop}%
  \fi \next}

\multibold
qwertyuiopasdfghjklzxcvbnmQWERTYUIOPASDFGHJKLZXCVBNM.

\multical
QWERTYUIOPASDFGHJKLZXCVBNM.

\multimathop
dist div dom meas sign Sign supp .


\def\accorpa #1#2{\eqref{#1}--\eqref{#2}}
\def\Accorpa #1#2 #3 {\gdef #1{\eqref{#2}--\eqref{#3}}%
  \wlog{}\wlog{\string #1 -> #2 - #3}\wlog{}}


\def\separa{\noalign{\allowbreak}}

\def\sign{\mathop{\rm sign}\nolimits}
\def\Sign{\mathop{\rm Sign}\nolimits}

\def\seto{\mathrel{{\scriptscriptstyle\searrow}}}

\def\graffe #1{\mathopen\{#1\mathclose\}}

\def\<#1>{\mathopen\langle #1\mathclose\rangle}
\def\norma #1{\mathopen \| #1\mathclose \|}

\def\iot {\int_0^t}
\def\ioth{\int_t^{t+h}\!\!}

\def\intQt{\int_{Q_t}}
\def\intQ{\int_Q}
\def\iO{\int_\Omega}

\def\dt{\partial_t}
\def\dn{\partial_n}

\def\cpto{\,\cdot\,}

\def\checkmmode #1{\relax\ifmmode\hbox{#1}\else{#1}\fi}
\def\aeO{\checkmmode{a.e.\ in~$\Omega$}}
\def\aeQ{\checkmmode{a.e.\ in~$Q$}}

\def\aeS{\checkmmode{a.e.\ on~$\Sigma$}}

\def\aat{\checkmmode{for a.a.~$t\in(0,T)$}}


\def\erre{{\mathbb{R}}}




\def\genspazio #1#2#3#4#5{#1^{#2}(#5,#4;#3)}
\def\spazio #1#2#3{\genspazio {#1}{#2}{#3}T0}

\def\L {\spazio L}
\def\H {\spazio H}
\def\W {\spazio W}

\def\C #1#2{C^{#1}([0,T];#2)}


\def\Lx #1{L^{#1}(\Omega)}
\def\Hx #1{H^{#1}(\Omega)}

\def\Cx #1{C^{#1}(\overline\Omega)}

\def\LQ #1{L^{#1}(Q)}

\def\CS #1{C^{#1}(\Sigma)}

\def\Luno{\Lx 1}
\def\Ldue{\Lx 2}
\def\Linfty{\Lx\infty}

\def\Huno{\Hx 1}
\def\Hdue{\Hx 2}


\def\LQ #1{L^{#1}(Q)}


\let\theta\vartheta
\let\eps\varepsilon
\let\phi\varphi

\let\TeXchi\chi                         
\newbox\chibox
\setbox0 \hbox{\mathsurround0pt $\TeXchi$}
\setbox\chibox \hbox{\raise\dp0 \box 0 }
\def\chi{\copy\chibox}


\def\QED{\hfill $\square$}


\def\thetaz{\theta_0}
\def\phiz{\phi_0}
\def\chiz{\chi_0}
\def\phistar{\phi^*}
\def\etastar{\eta^*}
\def\rhostar{\rho^*}
\def\xistar{\xi^*}

\def\Omegaz{\Omega_0}

\def\normaV #1{\norma{#1}_V}
\def\normaH #1{\norma{#1}_H}
\def\normaW #1{\norma{#1}_W}

\def\normaHeps #1{\norma{#1}_{H,\,\eps}}

\let\hat\widehat
\def\Beta{\hat{\vphantom t\smash\beta\mskip2mu}\mskip-1mu}
\def\Pi{\hat\pi}
\def\betaz{\beta^\circ}
\def\betaeps{\beta_\eps}
\def\signeps{\sign_\eps}
\def\Signeps{\Sign_\eps}

\def\thetaeps{\theta_\eps}
\def\phieps{\phi_\eps}
\def\xieps{\xi_\eps}
\def\etaeps{\eta_\eps}
\def\chieps{\chi_\eps}
\def\sigmaeps{\sigma_\eps}
\def\geps{g_\eps}
\def\psieps{\psi_\eps}

\def\az{a_0}
\def\bz{b_0}
\def\psiz{\psi_0}
\def\Tstar{T^*}



\def\CO{C_\Omega}
\def\CS{C_{str}}
\def\hatC{\hat C}

\def\Mz{M_0}
\def\Mpi{M_\pi^*}
\def\Mrho{M(\rho)}
\def\Arho{A(\rho)}
\def\Ctilde{\widetilde C}

\Begin{document}


\title{Sliding mode \betti{control} \\
for a \betti{nonlinear} phase-field system
}

\author{}
\date{}
\maketitle
\Bcenter
\vskip-2cm
{\large\sc Viorel Barbu$^{(1)}$}\\
{\normalsize e-mail: {\tt vb41@uaic.ro}}\\[.25cm]
{\large\sc Pierluigi Colli$^{(2)}$}\\
{\normalsize e-mail: {\tt pierluigi.colli@unipv.it}}\\[.25cm]
{\large\sc Gianni Gilardi$^{(2)}$}\\
{\normalsize e-mail: {\tt gianni.gilardi@unipv.it}}\\[.25cm]
{\large\sc Gabriela Marinoschi$^{(3)}$}\\
{\normalsize e-mail: {\tt gabriela.marinoschi@acad.ro}}\\[.25cm]
{\large\sc Elisabetta Rocca$^{(2)}$}\\
{\normalsize e-mail: {\tt elisabetta.rocca@unipv.it}}\\[.45cm]
$^{(1)}$
{\small ``Al. I. Cuza'' University, Ia\c si, Romania\\
and \ Romanian Academy, Ia\c si Branch\\
\gianniold{Bdul.\ Carol I 11, Ia\c si, Romania}}\\[.2cm]
$^{(2)}$
{\small Dipartimento di Matematica ``F. Casorati'', Universit\`a di Pavia}\\
{\small \gianniold{and \ IMATI-CNR, Pavia}}\\
{\small via Ferrata 1, 27100 Pavia, Italy}\\[.2cm]
$^{(3)}$
{\small ``Gheorghe Mihoc-Caius Iacob'' Institute of Mathematical Statistics\\
and \ Applied Mathematics of the Romanian Academy}\\
{\small Calea 13 Septembrie 13, 050711 Bucharest, Romania}\\[.2cm]
\Ecenter

\Begin{abstract}
\gianniold{In the present contribution the sliding mode control (SMC) problem for a phase-field model of Caginalp type is considered. 
First we prove the well-posedness and some regularity results for the phase-field type state systems 
modified by the state-feedback control laws. 
Then, we show that the chosen SMC laws force the system to reach within finite time the sliding manifold
(that we chose in order that one of the physical variables or a combination of them remains constant in time). 
We study three different types of feedback control laws: 
the first one appears in the internal energy balance 
and forces a linear combination of the temperature and the phase to reach a given (space dependent) value,
while the second and third ones are added in the phase relation and lead the phase  
onto a prescribed \revdue{target.}
While the control law is non-local in space for the first two problems,
it is local in the third one, i.e.,
its value at any point and any time just depends on the value of the state.}
\vskip3mm

\End{abstract}

\salta

\pagestyle{myheadings}
\newcommand\testopari{\sc Barbu \ --- \ Colli \ --- \ Gilardi \ --- \ Marinoschi \ --- \ Rocca}
\newcommand\testodispari{\sc Sliding modes for a phase-field system}
\markboth{\testodispari}{\testopari}

\finqui


\section{Introduction}
\label{Intro}
\setcounter{equation}{0}

\betti{%
Sliding mode control (SMC) has for many years been \recogniz ed as one of the
fundamental approaches for the systematic design of robust controllers for nonlinear complex
dynamic systems that operate under uncertainty. 
Moreover, SMC is  nowadays considered a classical tool for the regulation 
of continuous - or discrete - time systems in finite-dimensional settings 
(cf., e.g., the monographs \cite{BFPU08, EFF06, ES99, FMI11, I76, Utkin92, UGS09, YO99}).}

\betti{%
The main advantage of sliding mode control is that it  allows the separation of the motion of the overall system
in independent partial components of lower dimensions, and consequently it reduces the
complexity of the control problem. 
The design of feedback control systems with sliding modes implies the design of suitable control functions enforcing motions along ad-hoc manifolds. 
Hence, the main idea behind this scheme is first to identify a manifold of lower dimension 
(called the sliding manifold) 
where the control goal is fulfilled and such that the original system restricted to this sliding manifold 
has a desired \bhv, 
and then to act on the system through the control in order to constrain the evolution on it, 
that is\gianniold, to design a SMC-law that forces} 
\gianniold{the trajectories of the system to reach the sliding surface and maintains them on~it.}

\betti{%
Sliding mode controls, while being relatively easy to design, feature properties 
of both robustness with respect to unmodelled dynamics and insensitivity to external disturbances 
that are quite attractive in many applications. 
Hence, in the last years there has been a growing interest in the extension of the well developed methods 
for finite-dimensional systems described by ODEs 
(cf., e.g., \cite{LO02, O83, O00, OU83}) 
to the control of infinite-dimensional dynamical systems. 
While in some early works \cite{OU83, OU87, OU98} 
only special classes of evolutions were considered, 
the theoretical development in a general Hilbert space setting or for PDE systems has gained attention only in the last ten years. 
In this respect, we can quote the papers \cite{CRS11}, \cite{Levaggi13}, and \cite{XLGK13}  
dealing with sliding modes control for semilinear PDEs.  
In particular, in \cite{CRS11} the stabilization problem of a one-dimensional unstable heat conduction system
(rod) modeled by a parabolic partial differential equation, powered with a Dirichlet type actuator
from one of the boundaries was considered. 
A~delay-independent SMC strategy was proposed in \cite{XLGK13} to control a class of
quasi-linear parabolic PDE systems with time-varying delay, 
while in \cite{Levaggi13} the authors study a sliding mode control law for
a class of parabolic systems where the control acts through a Neumann boundary condition and the control space is finite-dimensional.}

\betti{%
In the present contribution we would like to employ -- to the best of our knowledge for the first time in the literature -- 
a~SMC technique \pier{for} a nonlinear PDE system of phase-field type. In particular, we consider the following rather simple version of the phase-field system of Caginalp type (see~\cite{Cag}):}
\Bsist
  & \dt \bigl( \theta + \ell \phi \bigr) - \kappa \Delta\theta = f
  & \quad \hbox{in $Q:=(0,T)\times\Omega$}
  \label{Iprima}
  \\
  & \dt\phi - \nu \Delta\phi + F'(\phi) = \gamma \theta 
  & \quad \hbox{in $Q$}
  \label{Iseconda}
\Esist
where $\Omega$ is the \rev{three-dimensional} domain \pier{in which} the evolution takes place,
$T$~is some final time,
$\theta$~denotes the relative temperature around some critical value
that is taken to be $0$ without loss of generality,
and $\phi$ is the order parameter.
Moreover, $\ell$, $\kappa$, $\nu$ and $\gamma$ are positive constants,
$f$~is a source term
and $F'$ \rev{represents} the derivative of a double-well potential~$F$.
Typical examples are
\begin{align}
  & F_{reg}(r) = \gianniold{\frac 14 \, (r^2-1)^2} \,,
  \quad r \in \erre
  \label{regpot}
  \\
  & F_{log}(r) = \bigl( (1+r)\ln (1+r)+(1-r)\ln (1-r) \bigr) - c_0 \, r^2 \,,
  \quad r \in (-1,1)
  \label{logpot}
  \\
  & F_{obs}(r) = I(r) - c_0 \, r^2 \,,
  \quad r \in \erre
  \label{obspot}
\end{align}
where $c_0>1$ in \eqref{logpot} in order to produce a double well,
while $c_0$ is an arbitrary positive number in~\eqref{obspot},
and the function $I$ in~\eqref{obspot} is the indicator function of~$[-1,1]$, i.e.,
it takes the values $0$ or $+\infty$ according to whether or not $r$ belongs to~$[-1,1]$.
The potential \eqref{regpot} and \eqref{logpot} are 
\gianniold{the usual classical regular potential} and the so-called logarithmic potential, respectively.
More generally, the potential $F$ could be just the sum 
$$\rev{F=\Beta+\Pi,}$$
where $\Beta$ is a convex function that is allowed to take the value~$+\infty$,
and $\Pi$ is a smooth perturbation (not necessarily concave).
In such a case, $\Beta$~is supposed to be proper and lower semicontinuous
so that its subdifferential is well-defined and can replace the derivative
which might not exist.
This happens in the case~\eqref{obspot}
and equation \eqref{Iseconda} becomes a differential inclusion.

The above system is complemented by initial conditions like $\theta(0)=\thetaz$ and $\phi(0)=\phiz$
and suitable boundary conditions.
\pier{Concerning the latter, as very usual}
we take the homogeneous Neumann condition for both $\theta$ and~$\phi$, that~is,
\Beq
  \dn\theta = 0
  \aand
  \dn\phi = 0
  \quad \hbox{on $\Sigma := (0,T)\times\Gamma $}
  \non
\Eeq
where $\Gamma$ is the boundary of~$\Omega$
and $\dn$ is the (say, outward) normal derivative.

\pier{%
Equations \eqref{Iprima}--\eqref{Iseconda} yield a system of phase field type. Such systems 
have been introduced (cf.~\cite{Cag}) in order to include phase dissipation 
effects in the dynamics of moving interfaces arising in thermally induced 
phase transitions. In our case, we move from the following expression for the 
total free energy
\begin{equation}
\calF (\theta, \varphi) = \int_\Omega \left( - \frac{c_0}2 \theta^2 - \gamma\theta \varphi + F (\varphi) + \frac\nu2 |\nabla \varphi |^2 \right)
\label{free}
\end{equation}
where $c_0 $ and $\gamma$ stand for specific heat and latent heat coefficients,
respectively, with a terminology motivated by earlier studies (see~\cite{duvaut})
on the Stefan problem; we refer to the monography \cite{fremond} 
which deals with phase change models as well.  In this connection, let us
introduce the enthalpy $e$   by
$$ 
e= - \frac{\delta \calF}{\delta \theta} \quad (- \hbox{ the variational derivative of $\calF $ with respect to } \theta )
$$    
that is $e=c_0 \theta +\gamma \varphi$. Then, the governing balance and phase equations 
are given~by 
\begin{gather}
\dt e + \div {\bf q}  = \tilde f 
\label{1phys}
\\
\dt \varphi +  \frac{\delta \calF}{\delta\varphi} =0
\label{2phys}
\end{gather}
where ${\bf q} $ denotes the thermal flux vector,  $\tilde f$ represents some heat source and  the variational derivative of $\calF$ with  respect to $\varphi$ appears in 
\eqref{2phys}. Hence, \eqref{2phys} reduces 
exactly to \eqref{Iseconda} along with the homogeneous Neumann boundary condition for $\varphi$. 
Moreover, if we assume the classical Fourier law $ {\bf q} = - \tilde\kappa \, \nabla \theta $,
then  \eqref{1phys} is nothing but the usual energy balance equation of the 
Caginalp model~\cite{Cag}. By setting $\ell :=\gamma/c_0$,  $\kappa:= \tilde\kappa /c_0$,
$f := \tilde f /c_0$,  we easily see that \eqref{Iprima} follows from  \eqref{1phys} and 
the Neumann boundary condition for $\theta$ is a consequence of the no-flux condition $ {\bf q} \cdot {\bf n} =0 $ on the boundary.  We also point out that the above phase field system has received a good deal of attention in the last decades and it 
can be deduced as a special gradient-flow problem (cf., e.g., \cite{RS06} and references therein).
}%

\betti{As \pier{already noticed,} the well-posedness, the long-time \pier\bhv\ of solutions, 
and also the related optimal control problems have been widely \rev{studied} in the literature. 
We refer, without \pier{any} sake of completeness, e.g., to \cite{BrokSpr, EllZheng, GraPetSch, KenmNiez, Lau} 
and references therein for the well-posedness \pier{and long time \bhv}\ results and to 
\pier{\cite{CGM, CoGiMaRo, HoffJiang, HKKY}}
for the related optimal control problems.
}%

The present paper is also related to the control problems, 
but it goes in the direction of designing sliding mode controls for the above phase-field system. 
Indeed our main objective is to find out some state-feedback control laws 
$(\theta,\phi)\mapsto u(\theta,\phi)$ that can be inserted in one of the equations
in order that the dynamics of the system modified in this way
forces the value $(\theta(t),\phi(t))$ of the solution 
to reach some manifold of the phase space in a finite time
and then lie there with a sliding mode. 

The first analytical difficulty consists in deriving the equations
governing the sliding modes and the conditions for this motion to exist. 
The problem needs the development of special methods, 
since the conventional theorems regarding existence and uniqueness 
of solutions  are not directly applicable. 
Moreover, we need to manipulate the system through the control in order 
to constrain the evolution on the desired sliding manifold.  
In particular, we study three cases. 

In the first one, 
a feedback control is added to the internal energy balance equation \eqref{Iprima} 
in order to force a linear relationship between $\theta$ and~$\phi$; 
in the second case,
a prescribed \gianniold{distribution} $\phistar$ of the order parameter is forced by means of a feedback control  
added to the phase dynamics~\eqref{Iseconda}. 
Notice that both these choices can be considered physically meaningful 
in the framework of phase transition processes, 
since in both cases the quantities we are forcing to reach time-independent values
may have a physical meaning.
In the first problem, we can take the internal energy as a particular case,
while the target $\phistar$ we force for the phase parameter in the second problem
could represent one of the so called pure phases 
(e.g., pure water or pure ice in a water-ice phase change process). 
Moreover, in both cases we have reduced the problem to a simplified \pier{dynamics} involving only 
the evolution of $\varphi$ in the first case and only of $\theta$ in the second one
(cf.\ also Remark~\ref{Simpldyn}).

In each of the above problems, the control law we introduce is non-local in space, i.e.,
the value at $(t,x)$ of the control depends on the whole state $(\theta(t,\cpto),\phi(t,\cpto))$
at the time~$t$ and not only on the value $(\theta(t,x),\phi(t,x))$.
The objective of the third problem 
is to design a control law that reaches the same target as in the second one and is local at the same time.
However, such a problem looks much more difficult
and we can ensure the existence of the desired sliding mode only under 
\gianni{a suitable compatibility condition on~$\Omega$.}

The paper is \organiz ed as follows.
In the next section, we list our assumptions, state the problem in a precise form
and present our results.
The last two sections are devoted to the corresponding proofs.
Section~\ref{WELLPOSEDNESS} deals with \gianniold{well-posedness and regularity},
while the existence of the sliding modes is proved in Section~\ref{SLIDINGMODES}.


\section{Statement of the problem and results}
\label{STATEMENT}
\setcounter{equation}{0}

In this section, we describe the problem under study
and present our results.
As in the Introduction,
$\Omega$~is the body where the evolution takes place.
We assume 
$$\rev{\Omega\subset\erre^3 \ \hbox{ 
to~be open, bounded, connected, and smooth}}$$
and \rev{write} $|\Omega|$ for its Lebesgue measure.
Moreover, $\Gamma$ and $\dn$ still stand for
the boundary of~$\Omega$ and the outward normal derivative, respectively.
Given a finite final time~$T>0$,
we set for convenience
\gianniold{$Q:=(0,T)\times\Omega$}.
Now, we specify the assumptions on the structure of our system.
We assume that
\begin{align}
  & \ell, \, \kappa, \, \nu, \, \gamma \in (0,+\infty)
  \label{hpconst}
  \\
  & \Beta : \erre \to [0,+\infty]
  \quad \hbox{is convex, proper and l.s.c.}
  \quad \hbox{with} \quad
  \Beta(0) = 0
  \label{hpBeta}
  \\
  & \Pi: \erre \to \erre
  \quad \hbox{is a $C^1$ function}
  \aand
  {\Pi\,}'
  \quad \hbox{is \rev{uniformly \Lip}.}
  \label{hpPi}
\end{align}
\Accorpa\HPstruttura hpconst hpPi
We set for brevity
\Beq
  \beta := \partial\Beta 
  \aand
  \pi := {\Pi\,}'
  \label{defbetapi}
\Eeq
and denote by $D(\beta)$ and $D(\Beta)$
the effective domains of $\beta$ and~$\Beta$, respectively.
Next, in order to simplify notations, we~set
\Beq
  V := \Huno, \quad
  H := \Ldue, \quad
  W := \graffe{v\in\Hx2: \dn v=0}
  \label{defspazi}
\Eeq
\gianniold{and endow the spaces $V$ and $H$ with their standard norms
$\normaV\cpto$ and~$\normaH\cpto$.
On the contrary, we write $\normaW\cpto$ for the norm in $W$ defined~by
\Beq
  \normaW v^2 = \normaH v^2 + |\Omega|^{4/3} \normaH{\Delta v}^2
  \quad \hbox{for every $v\in W$}
  \label{defnormaW}
\Eeq
and we term $\CO$ the best constant realizing the inequality
\Beq
  \norma v_\infty \leq \CO \, \normaW v
  \quad \hbox{for every $v\in W$}.
  \label{embedding}
\Eeq
The reason of this choice will be explained later~on
(see the forthcoming Remark~\ref{Smallness}).
Now, we just notice that $\normaW\cpto$ is equivalent
to the norm induced on $W$ by the standard one in~$\Hdue$
(thanks to the regularity theory of elliptic equations)
and~that the constant $\CO$ actually exists 
due to the continuous embedding $\Hdue\subset\Cx0$
(since $\Omega\subset\erre^3$ is bounded and smooth)
and~only depends on~$\Omega$ \rev{(see, e.g., \cite{GT}).}
Finally, for the norms both in $\Lx\infty$ and in $\LQ\infty$
we use the same symbol $\norma\cpto_\infty$ whenever no confusion can arise.}

Furthermore,  \rev{\revdue{let}%
\Beq
  \Sign : H \to 2^H
  \quad \hbox{\revdue{be} the subdifferential of the map} \quad
  \normaH\cpto: H \to \erre
  \label{defSign}
\Eeq
i.e., 
\Bsist
  && \Sign v = \frac v {\normaH v}
  \quad \hbox{if $v\in H$ and $v\not=0$}
  \label{signv}
  \\
  && 
  \Sign 0 \enskip
  \hbox{is the closed unit ball of $H$.}
  \label{signz}
\Esist
}%
Thus, $\beta$ and $\Sign$ are maximal monotone operators
on $\erre$ and~$H$, respectively
(see, e.g., \cite[Thm.~2.8, p.~47]{Barbu}).
In the sequel, we use the same symbol $\beta$ to denote
the maximal monotone operator induced on $L^2$-spaces.

\step \rev{Yosida regularizations of $\beta$ and $\Sign$}

\rev{Let us} introduce the Yosida \regulariz ation 
$\betaeps:\erre\to\erre$ and $\Signeps:H\to H$ at level~\pier{$\eps>0$}
(see, e.g., \cite[formulas~(2.26), p.~37]{Barbu})
as~well as the Moreau \regulariz ation of~$\normaH\cpto$
(see, e.g., \cite[formula~(2.38), p.~48]{Barbu})
\Beq
  \normaHeps v := \min_{w\in H} \bigl\{ 
                   {\textstyle\frac 1 {2\eps}} \normaH{w-v}^2 + \normaH w
                 \bigr\}
  = \int_0^{\normaH v} \!\! \min \{ s/\eps , 1 \} \, ds                 
  \quad \hbox{for $v\in H$}.
  \label{moreau}
\Eeq
\pier{%
For the reader's convenience, we sketch the justification of the last equality of~\eqref{moreau}.
We write $\norma\cpto$ instead of $\normaH\cpto$ for simplicity.
For $w\in H$ and $y\geq0$ we~set
\Beq
  G(w) := \frac 1 {2\eps} \, \rev{\norma{w-v}}^2 + \rev{\norma w}
  \aand
  g(y) := \frac 1 {2\eps} (y-\norma v)^2 + y
  \non
\Eeq
and observe that the triangle inequality $\bigl|\norma w-\norma v\bigr|\leq\norma{w-v}$ yields
$G(w)\geq g(\norma w)$ for every $w\in H$.
Now, from one side, one easily checks that
\Beq
  \min_{y\geq 0} g(y) = \frac 1 {2\eps} \, \norma v^2
  \quad \hbox{if $\norma v\leq\eps$}
  \aand
  \min_{y\geq 0} g(y) = \norma v - \frac \eps 2
  \quad \hbox{if $\norma v>\eps$} .
  \non
\Eeq
This means that $\min_{y\geq 0} g(y)$ coincides with the \rhs\ of \eqref{moreau}.
On the other hand, we have
\Beq
  G(0) = \frac 1 {2\eps} \, \norma v^2
  \quad \hbox{in any case,}
  \aand
  G \bigl( (1-\eps/\norma v) v \bigr)
  = \norma v - \frac \eps 2
  \quad \hbox{if $\norma v>\eps$}.
  \non
\Eeq
Thus, $\min_{w\in H}G(w)=\min_{y\geq 0}g(y)$ and \eqref{moreau} is proved.
Next, we recall that
}%
$\betaeps$ and $\Signeps$ are monotone and that
(see, e.g., \cite[Prop.\ 2.2~(ii), p.~38]{Barbu} 
and \cite[Thm.~2.9, p.~48]{Barbu} for some of these properties) 
\Bsist
  && \hbox{$\Signeps v$\enskip is the gradient at $v$ of the $C^1$ functional\enskip $\normaHeps\cpto$}
  \label{dmoreau}
  \\
  && \Signeps v = \frac v {\max\{\eps,\normaH v\}}
  \quad \hbox{for every $v\in H$}
  \label{formulaSigneps}
  \\
  && \rev{(\Signeps v, v )_H \geq  
  \normaH{v} -\frac \eps 4 
  \quad \hbox{for every $v\in H$}} 
  \label{disG}
  \\
  && |\betaeps(r)| \leq |\betaz(r)|
  \quad \hbox{for every $r\in D(\beta)$}\gianniold, \quad \hbox{where}
  \non
  \\
  && \quad \betaz(r)\ \hbox{is the element of $\beta(r)$ having minimum modulus}\rev{.}
  \label{defbetaz}
\Esist
\rev{We point out that the Young inequality has been used to derive \eqref{disG}.}

At this point, we describe the state system 
modified by the state-feedback control law
and we study two cases.
In the first one, a~feedback control is added to the first equation~\eqref{Iprima}
in order to force a linear relationship between $\theta$ and~$\phi$;
in the second case, a~prescribed \gianniold{distribution} of the order parameter is forced
by means of a feedback control that is added to equation~\eqref{Iseconda}.
In principle, for the data, we require~that
\Bsist
  && f \in \LQ2 , \quad
  \thetaz \in V , \quad
  \phiz \in V
  \aand
  \Beta(\phiz) \in \Luno .
  \label{hpdati}
\Esist
\gianniold{Given} $\rho>0$ and some target that \pier{depends} on the case we want to consider,
we look for a quadruplet $(\theta,\phi,\xi,\sigma)$ satisfying
at least the regularity requirements
\Bsist
  && \theta, \phi \in \H1H \cap \L\infty V \cap \L2W
  \label{regsol}
  \\
  && \xi \in \L2H
  \aand
  \sigma \in \L\infty H ,
  \label{regxisigma}
\Esist
\Accorpa\Regsoluz regsol regxisigma
and solving the related system we introduce at once.
We notice that the homogeneous \gianniold{Neumann} boundary conditions for both $\theta$ and $\phi$
are contained in~\eqref{regsol} (see the definition \eqref{defspazi} of~$W$).
The problems corresponding to the cases sketched above are the following.

Given $\etastar\in W$ and $\alpha\in\erre$, the first system~is
\Bsist
  && \dt \bigl( \theta + \ell \phi \bigr) - \kappa \Delta\theta
  = f - \rho \sigma 
  \quad \aeQ
  \label{primaA}
  \\
  && \dt\phi - \nu \Delta\phi + \xi + \pi(\phi)
  = \gamma \theta
  \quad \aeQ
  \label{secondaA}
  \\
  && \xi \in \beta(\phi) \quad \aeQ
  \label{terzaA}
  \\
  && \sigma(t) \in \Sign(\theta(t) + \alpha\phi(t) - \etastar)
  \quad \aat 
  \label{quartaA}
  \\
  && \theta(0) = \thetaz
  \aand
  \phi(0) = \phiz \,.
  \label{cauchyA}
\Esist
\Accorpa\PblA primaA cauchyA
\gianniold{In the sequel, we also term such problem
{\sl Problem~(A)\/}.}

The second problem, which we call {\sl Problem~(B)}\/, 
depends on a given \gianniold{$\phistar\in W$ and} consists in the equations
\Bsist
  && \dt \bigl( \theta + \ell \phi \bigr) - \kappa \Delta\theta 
  =  f
  \quad \aeQ
  \label{primaB}
  \\
  && \dt\phi - \nu \Delta\phi + \xi + \pi(\phi)
  = \gamma \theta - \rho \sigma 
  \quad \aeQ
  \label{secondaB}
  \\
  && \xi \in \beta(\phi) \quad \aeQ
  \label{terzaB}
  \\
  && \sigma(t) \in \Sign(\phi(t) - \phistar)
  \quad \aat 
  \label{quartaB}
  \\
  && \theta(0) = \thetaz
  \aand
  \phi(0) = \phiz \,.
  \label{cauchyB}
\Esist
\Accorpa\PblB primaB cauchyB

\gianniold{The last case, \pier{termed} {\sl Problem~(C)}\/,
is the same as the previous one with the following difference:
the non-local operator $\Sign$ is replaced by the local $\sign:\erre\to 2^{\erre}$
defined~by
\Beq
  \sign r := \frac r {|r|}
  \quad \hbox{if $r\not=0$}
  \aand
  \sign 0 := [-1,1] .
  \label{defsign}
\Eeq
Notice that $\sign$ is the subdifferential of the real function
$r\mapsto|r|$ and thus is maximal monotone.
For the sake of clarity, we write Problem~(C), \pier{explicitly}.
Given~\gianni{$\phistar\in W$}, 
we look for a quadruplet $(\theta,\phi,\xi,\sigma)$ satisfying
\Bsist
  && \dt \bigl( \theta + \ell \phi \bigr) - \kappa \Delta\theta 
  =  f
  \quad \aeQ
  \label{primaC}
  \\
  && \dt\phi - \nu \Delta\phi + \xi + \pi(\phi)
  = \gamma \theta - \rho \sigma 
  \quad \aeQ
  \label{secondaC}
  \\
  && \xi \in \beta(\phi) \quad \aeQ
  \label{terzaC}
  \\
  && \pier{\sigma \in \sign(\phi - \phistar)
  \quad \aeQ}
  \label{quartaC}
  \\
  && \theta(0) = \thetaz
  \aand
  \phi(0) = \phiz \,.
  \label{cauchyC}
\Esist
\Accorpa\PblC primaC cauchyC
Here are our results on the well-posedness of the above problems.
}%

\Bthm
\label{WellposednessA}
Assume \HPstruttura, \eqref{hpdati}\gianniold,
\Beq
  \etastar \in W
  \aand
  \alpha \in \erre .
  \label{hpexA}
\Eeq
Then, for every $\rho>0$, \pier{Problem~(A)} has at least a solution $(\theta,\phi,\xi,\sigma)$
satisfying \Regsoluz\ and the estimates
\Bsist
  && \norma\theta_{\L\infty H\cap\L2V}
  + \norma\phi_{\H1H\cap\L2W}
  \non
  \\
  && \quad {}
  + \norma\xi_{\L2H}
  + \betti{\norma\sigma_{\L\infty H}}
  \leq C_1
  \qquad
  \label{stimaAuno}
  \\
  && \norma\theta_{\H1H\cap\L2W}
  \leq C_2 \bigl( \rho^{1/2} + 1 \bigr)
  \label{stimaAdue}
\Esist
where $C_1$ and $C_2$ depend only on the quantities involved
in assumptions \HPstruttura, \eqref{hpdati} and~\eqref{hpexA}.
Moreover, the solution is unique if $\alpha=\ell$.
Furthermore, if in addition
\Beq
  \phiz \in W
  \aand
  \betaz(\phiz) \in H 
  \label{hppiuregA}
\Eeq
then, there exists a solution that also satisfies
\Bsist
  && \phi \in \W{1,\infty}H \cap \H1V \cap \L\infty W
  \aand
  \xi \in \L\infty H
  \qquad
  \label{phipiureg}
  \\
  && \norma\phi_{\W{1,\infty}H\cap\H1V\cap\L\infty W}
  \leq C_3 \bigl( \rho^{1/2} + 1 \bigr)
  \label{stimaApiureg} 
\Esist
where $C_3$ depends on the norms related to~\eqref{hppiuregA} as well.
\gianniold{In particular, $\phi$~is bounded.
Finally, the component $\theta$ of any solution 
satisfying all the above regularity requirements
is bounded whenever $\thetaz\in V\cap\Linfty$ and $f\in\L\infty H$.}
\Ethm

\Bthm
\label{WellposednessB}
Assume \HPstruttura, \eqref{hpdati}\gianniold{, as well as}
\Beq
  \gianniold{\phistar\in W \aand \betaz(\phistar) \in H} \,.
  \label{hpexB}
\Eeq
Then, for every $\rho>0$, \pier{Problem~(B)} has at least a solution $(\theta,\phi,\xi,\sigma)$
satisfying \Regsoluz\ and the estimates
\Bsist
  && \norma\theta_{\L\infty H\cap\L2V}
  + \norma\phi_{\L\infty H\cap\L2V}
  + \betti{\norma\sigma_{\L\infty H}}
  \leq C_4
  \label{stimaBuno}
  \\
  && \norma\theta_{\H1H\cap\L2W}
  + \norma\phi_{\H1H\cap\L2W}
  \non
  \\
  && \quad {}
  + \gianniold{\norma{\xi+\rho\sigma}_{\L2H}}
  \leq C_5 \bigl( \rho^{1/2} + 1 \bigr) 
  \label{stimaBdue}
\Esist
where $C_4$ and $C_5$ depend only on the quantities involved
in assumptions \HPstruttura, \eqref{hpdati} and~\eqref{hpexB}.
Furthermore, the components $\theta$ and $\gianniold\phi$ of the solution
are uniquely determined,
and $\xi$ and $\sigma$ are uniquely determined as well if $\beta$ is single-valued.
\Ethm

\gianniold{A similar result holds for Problem~(C).
We present the corresponding statement in a more accurate form
for a reason that will be clear later~on.}

\gianniold{%
\Bthm
\label{WellposednessC}
Assume \HPstruttura, \eqref{hpdati} and~\eqref{hpexB}.
Then, for every $\rho>0$, \pier{Problem~(C)} has at least a solution $(\theta,\phi,\xi,\sigma)$
satisfying \Regsoluz.
Furthermore, the components $\theta$ and $\phi$ of the solution
are uniquely determined,
and $\xi$ and $\sigma$ are uniquely determined as well if $\beta$ is single-valued.
Finally, if the conditions
\Beq
  f \in \H1H , \quad
  \thetaz \in W , \quad
  \phiz \in W
  \aand
  \betaz(\phiz) \in H 
  \label{hppiuregC}
\Eeq
are assumed in addition,
then \eqref{phipiureg} holds as well as
\Beq
  \theta \in \W{1,\infty}H \cap \H1V \cap \L\infty W .
  \label{thetapiureg}
\Eeq
In particular, both $\phi$ and $\theta$ are bounded.
Moreover, the estimates
\Bsist
  && \norma{\phi-\phistar}_\infty
  \leq \rho \, \CS \, \CO |\Omega|^{7/6}
  + C_6 
  \label{chibddC}
  \\
  && \norma\theta_\infty 
  \leq \rho \, \CS \, \CO |\Omega|^{7/6}
  + C_7
  \quad
  \label{thetabddC}
\Esist
hold true with a structural constant $\CS$ depending only 
on the physical parameters $\ell$, $\kappa$, $\nu$ and~$\gamma$,
the constant $\CO$ given by~\eqref{embedding}
and some constants $C_6$ and $C_7$ depending 
on the structure of the systems, $\Omega$, $T$ and on the norms of the data involved.
\Ethm
}%

\Brem
The above results are quite general.
In particular, both potentials \eqref{regpot} and \eqref{logpot} 
are certainly allowed
and the multi-valued potential \eqref{obspot} has to be excluded
just in the parts of Theorems~\ref{WellposednessB} and~\ref{WellposednessC}
regarding uniqueness for the pair $(\xi,\sigma)$,
which might be not uniquely determined, in general.
Concerning the constant~$\CS$ of Theorem~\ref{WellposednessC},
we will prove that we can take
\Beq
  \CS = 2 \max \Bigl\{ \frac {6^{1/2}} {\gianni\nu} \,,\, \frac \ell {\kappa^{1/2} \nu^{1/2}} + \frac {4\ell} \kappa \Bigr\}.
  \label{costC}
\Eeq
However, no sharpness is guaranteed at all.
\Erem

For each of the first two problems, the existence of the desired sliding mode is ensured for $\rho$ large enough.
\gianni{For every $T>0$ we have indeed}

\Bthm
\label{SlidingA}
Assume \HPstruttura, \eqref{hpdati}, \eqref{hpexA}, \gianniold{\eqref{hppiuregA} and $f\in\L\infty H$}.
Then, for some $\rhostar>0$ and for every $\rho>\rhostar$,
there exist a solution $(\theta,\phi,\xi,\sigma)$
to problem \PblA\ and a time $\Tstar\in\gianniold[0,T)$ such~that
\Beq
  \theta(t) + \alpha \phi(t) = \etastar
  \quad \aeO
  \quad \hbox{for every $t\in[\Tstar,T]$}.
  \label{slidingA}
\Eeq
\Ethm

\Bthm
\label{SlidingB}
Assume \HPstruttura, \eqref{hpdati} and \eqref{hpexB}.
Then, for some $\rhostar>0$ and for every $\rho>\rhostar$,
there exist a solution $(\theta,\phi,\xi,\sigma)$
to problem \PblB\ and a time $\Tstar\in\gianniold[0,T)$ such~that
\Beq
  \phi(t) = \phistar
  \quad \aeO
  \quad \hbox{for every $t\in[\Tstar,T]$}.
  \label{slidingB}
\Eeq
\Ethm

\Brem
\label{CalcolostarAB}
In the proof we give in Section~\ref{SLIDINGMODES},
we compute possible values of $\rhostar$ and $\Tstar$
that fit the conclusions of our results.
\gianni{For Problems~(A) and~(B), we can~take respectively}
\Bsist
  && \gianni\rhostar := C_A^2 + 2C_A + \frac 2T \, \normaH{\thetaz+\alpha\phiz-\etastar}
  \aand
  \gianni\Tstar := \frac {2 \normaH{\thetaz+\alpha\phiz-\etastar}} {\rho - C_A^2 - 2C_A}
  \non
  \\
  && \gianni\rhostar := 2C_B + \frac 2T \, \normaH{\phiz-\phistar}
  \aand
  \gianni\Tstar := \frac {2 \normaH{\phiz-\phistar}} {\rho - 2C_B}
  \non
\Esist
where the constants $C_A$ and~$C_B$ 
are constructed in the proofs of Theorems \ref{WellposednessA} and~\ref{WellposednessB}
in order that
\Bsist
  && \norma{f-(\ell-\alpha)\dt\phi-\kappa\alpha\Delta\phi-\kappa\Delta\etastar}_{\L\infty H}
  \leq C_A \bigl( \rho^{1/2} + 1 \bigr)
  \non
  \\
  && \norma{\gamma\theta+\nu\Delta\phistar-\betaz(\phistar)\rev{{}- \pi (\phi)} }_{\L\infty H}
  \leq C_B \,.
  \non
\Esist
\gianni{More precisely, we refer to \accorpa{primaAter}{defrhostarA} 
and \accorpa{secondaBter}{defrhostarB} and we notice that
our starting point in those proofs
is the validity of the \rev{analogous} estimates for the solutions to the approximating problems
obtained by replacing the monotone operators by their Yosida \regulariz ations.} \rev{It follows that the above values of $\rhostar$ and $\Tstar$ depend continuously on the potentials and on the physical parameters of the systems.}
\rev{We also observe} that the time $\Tstar$ is roughly proportional to~$1/\rho$ in both cases,
whence it tends to zero as $\rho$ tends to infinity\rev{, i.e.,} the sliding mode can be forced to start as soon as one desires
by prescribing a sufficiently big factor $\rho$ in front of the feedback control. 
\Erem

\Brem
\label{Simpldyn}
The minimal value of $\,\Tstar$ of the first statement 
(if~it is positive) also satisfies the following property:
the function $t\mapsto\normaH{\theta(t) + \alpha \phi(t) - \etastar}$ is strictly decreasing on~$[0,\Tstar]$.
\gianniold{A~similar remark holds for the function $t\mapsto\normaH{\phi(t)-\phistar}$ 
in the second statement 
(and in the next one, at least under some reinforcement of the assumptions,
as shown in Remark~\ref{DecreasingC})}.
In each case, the dynamics of the system is simpler after the time~$\,\Tstar$,
since one of the unknowns can be eliminated by using the sliding mode condition.
\gianniold{For instance}, in the second situation, the evolution of $\theta$ after $\Tstar$
is ruled just by the heat equation.
\Erem

The situation for Problem~(C) is different,
since we can ensure the existence of the desired sliding mode for $\rho$ large enough 
only if further conditions are fulfilled.
\gianni{Namely, we need a restriction
involving the structure of the system and the domain~$\Omega$
(that is why we have written the statement of Theorem~\ref{WellposednessC} in that form)}.
Our result only involves the component $\phi$ of the solution,
and we recall that~$\phi$ is uniquely determined.

\Bthm
\label{SlidingC}
Assume \HPstruttura, \eqref{hpdati}, \eqref{hpexB}, \eqref{hppiuregC} and
\Beq
  \Delta\phistar \in \Linfty
  \aand
  \betaz(\phistar) \in \Linfty .
  \label{hpslidingC}
\Eeq
Let $\CS$ and $\CO$ be the constants appearing 
in \eqref{thetabddC} and in~\eqref{embedding}, respectively,
and assume that
\Beq
  \gamma \, \CS\CO |\Omega|^{7/6} < 1 .
  \label{smallness}
\Eeq
Then, for some $\rhostar>0$ and for every $\rho>\rhostar$,
the following is true:
if $(\theta,\phi,\xi,\sigma)$ is a solution to problem \PblC,
there exists a time $\Tstar\in[0,T)$ such~that
\Beq
  \phi(t) = \phistar
  \quad \aeO
  \quad \hbox{for every $t\in[\Tstar,T]$}.
  \label{slidingC}
\Eeq
\Ethm

\Brem
\label{CalcolostarC}
\gianniold{%
Assume that the constants $\CS$, $\CO$ and $C_7$ realize
the inequalities \eqref{thetabddC} and~\eqref{smallness}
(i.e.,~in contrast with the situation of Remark~\ref{CalcolostarAB},
just such inequalities are required as a starting point).
Then, as shown in the proof we perform in the last section,
possible values of $\rhostar$ and $\Tstar$
that fit the conclusion of the above theorem are given~by
(here $L$ is the \Lip\ constant of~$\pi$)
\Bsist
  && \rhostar := \frac
  {\gamma C_7 + \nu \norma{\Delta\phistar}_\infty + \norma\xistar_\infty + \Mpi + \Mz/T}
  {1 - \gamma \, \CS\CO |\Omega|^{7/6}} 
  \aand
  \Tstar := \frac \Mz {\rho-\Arho} 
  \non
  \\
  \noalign{\smallskip}
  &&
  \hbox{where} \quad
  \Mpi := L(\Mz+\norma\phistar_\infty) + |\pi(0)| \,, \quad
  \Mz := \norma{\phiz-\phistar}_\infty 
  \aand
  \non
  \\
  && \Arho := \gamma \bigl( \CS\CO |\Omega|^{7/6} \rho + C_7  \bigr)
  + \nu \norma{\Delta\phistar}_\infty + \norma\xistar_\infty
  + \Mpi \,.
  \non
\Esist
In particular, the last \rev{two sentences} of Remark~\ref{CalcolostarAB} also \rev{apply} to the present case.}
\Erem

\Brem
\label{Smallness}
In order to understand the meaning of~\eqref{smallness},
let us assume that the structure of the system is chosen,
so that the physical constants are fixed,
and let us think of a class of open sets having the same shape.
Precisely, we fix an open set $\Omegaz$ of measure~$1$
and assume that
$\Omega=x_0+\lambda R\,\Omegaz$ for some $x_0\in\erre^3$, $\lambda>0$ and some rotation $R\in SO(3)$.
Then $\pier\lambda=|\Omega|^{1/3}$ and one easily checks that
our definition \eqref{defnormaW} of~$\normaW\cpto$ yields 
$\CO=C_{\Omegaz}\,|\Omega|^{-1/2}$,
since  the $H$-norms
of $v$ and of $\Delta v$ are properly balanced in the norm~$\normaW v$
under a rescaling of a function~$v$.
Then, the smallness condition \eqref{smallness}
means that $|\Omega|$ is small enough.
Indeed, the \lhs\ of \eqref{smallness}
\gianniold{becomes $\gamma\CS C_{\Omega_0}|\Omega|^{2/3}$}
in the chosen class of domains.
\Erem

In performing our a priori estimates in the remainder of the paper, 
we often account for the \Holder\ inequality and the elementary \gianniold{inequalities
(for arbitrary $a,b\geq0$)
\Beq
  (a+b)^{1/2}\leq a^{1/2}+b^{1/2} , \quad
  (a+b)^2 \leq 2a^2 + 2b^2
  \aand
  ab \leq \delta a^2 + \frac 1 {4\delta} \, b^2
  \label{elementary}
\Eeq
where $\delta>0$ in the latter (Young's inequality).
Moreover, we} repeatedly use the notation
\Beq
  Q_t := (0,t) \times \Omega \,.
  \label{defQ}
\Eeq
\gianni{For simplicity, we usually omit $dx$, $ds$, etc.\ in integrals.
More precisely, we explicitly write, e.g., $ds$ only if the variable $s$ actually appears 
in the function under the integral sign.}
Finally, while a particular care is taken in computing some constants,
we follow a general rule to denote less important ones,
in order to avoid boring calculations.
The small-case symbol $c$ stands for different constants \rev{independent of $\rho$ but depending
on~$\Omega$, the final time~$T$, the shape of the nonlinearities
and on the constants and the norms of
the functions involved in the assumptions of our statements.}
The dependence on $\rho$ will be always written \pier{explicitly}, indeed.
Hence, the meaning of $c$ might change from line to line 
and even in the same chain of equalities or inequalities.
On the contrary, we mark precise constants which we can refer~to
by using different symbols, e.g., capital letters,
\gianniold{mainly with indices, like in~\eqref{embedding}.}


\section{Proof of the well-posedness results}
\label{WELLPOSEDNESS}
\setcounter{equation}{0}

This section is devoted to the proof of \gianniold{Theorems~\ref{WellposednessA}--\ref{WellposednessC}}.
However, as far as existence is concerned, 
we confine ourselves to derive the formal a~priori estimates
that lead to the desired regularity
and just sketch how a completely rigorous proof could be performed.

\subsection{\rev{Proof of Theorem~\ref{WellposednessA}}}

We start with problem \PblA\ and transform it into an equivalent system in new unknown functions.
\rev{In order to argue in terms of the variable which the operator $\Sign$ applies to,} we~set
\Beq
  \eta := \theta + \alpha\phi - \etastar
  \label{defeta}
\Eeq
then, $\eta$ has to satisfy the \analogue\ of \eqref{regsol}
and the new problem is the following
\Bsist
  && \dt \bigl( \eta + (\ell-\alpha) \phi \bigr) - \kappa \Delta\eta + \kappa\alpha \Delta\phi
  = f \pier{{}+{}} \kappa \Delta\etastar - \rho \sigma 
  \quad \aeQ
  \label{primaAbis}
  \\
  && \dt\phi - \nu \Delta\phi + \xi + \pi(\phi)
  = \gamma (\eta - \alpha \phi + \etastar)
  \quad \aeQ
  \label{secondaAbis}
  \\
  && \xi \in \beta(\phi) \quad \aeQ
  \label{terzaAbis}
  \\
  && \sigma(t) \in \Sign(\eta(t))
  \quad \aat 
  \label{quartaAbis}
  \\
  && \eta(0) = \thetaz + \alpha\phiz - \etastar
  \aand
  \phi(0) = \phiz \,.
  \label{cauchyAbis}
\Esist
\Accorpa\PblAbis primaAbis cauchyAbis

\step
First a priori estimate

We multiply \eqref{primaAbis} and \eqref{secondaAbis} by $\eta$ and~$\dt\phi$, respectively,
sum up and integrate over $Q_t$ with an arbitrary $t\in(0,T]$.
Then, we add $\nu\intQt\phi\dt\phi=(\nu/2)\iO(|\phi(t)|^2-|\phiz|^2)$ to both sides.
\pier{With the help of \eqref{hpdati} and~\eqref{defSign}, we infer that}
\Bsist
  && \frac 12 \iO |\eta(t)|^2
  + \kappa \intQt |\nabla\eta|^2
  - \kappa\alpha \intQt \nabla\phi \cdot \nabla\eta
  + \rho \iot \normaH{\eta(s)} \, ds
  \non
  \\
  && \quad {}
  + \intQt |\dt\phi|^2
  + \rev{\frac \nu 2}\iO |\nabla\phi(t)|^2
  + \iO \Beta(\phi(t))
  + \frac \nu 2 \iO |\phi(t)|^2
  \non
  \\
  && \mathrel{\pier\leq} c
  - (\ell-\alpha) \intQt \dt\phi \, \eta
  + \intQt \bigl( f \pier{{}+{}} \kappa \Delta\etastar \bigr) \eta
  \non
  \\
  && \quad {}
  - \intQt \pi(\phi) \, \dt\phi
  + \gamma \intQt \bigl( \eta - \alpha \phi + \etastar \bigr) \dt\phi 
  + \nu \intQt \phi \dt\phi .
  \non
\Esist
Now, it is \sfw\ to use the linear growth of $\pi$ that follows from \Lip\ continuity, 
the Young and \Holder\ inequalities, \gianni{\eqref{hpdati}, \eqref{hpexA},}
and the Gronwall lemma to deduce~that
\Beq
  \norma\eta_{\L\infty H\cap\L2V}
  + \norma\phi_{\H1H\cap\L\infty V}
  + \norma{\Beta(\phi)}_{\L\infty\Luno}
  \leq c \,.
  \label{primastimaA}
\Eeq

\step 
Second a priori estimate

We write \eqref{secondaAbis} as
\Beq
  - \nu \Delta\phi(t) + \xi(t) = g_1(t)
  \aand
  \xi(t) \in \beta(\phi(t))
  \quad \aat
  \non
\Eeq
with an obvious meaning of $g_1$
and treat $t$ as a parameter.
We formally multiply by $\Delta\phi(t)$
(the~correct proof deals with the \regulariz ed problem)
and find $\normaH{\Delta\phi(t)}\leq\normaH{g_1(t)}$ \aat.
Then, we use \pier{\eqref{primastimaA}, \eqref{hpPi}, \eqref{hpexA}
(which imply $\norma{g_1}_{\L2H}\leq c$)},
elliptic regularity and a comparison in the above equation,
\pier{in order to} conclude that
\Beq
  \norma\phi_{\L2W} + \norma\xi_{\L2H} \leq c \,.
  \label{secondastimaA}
\Eeq

\step
Third a priori estimate

We write \eqref{primaAbis} as
\Beq
  \dt\eta - \kappa\Delta\eta + \rho\sigma = g_2
  \quad \hbox{with} \quad
  \norma{g_2}_{\L2H} \leq c
  \label{perterzaA}
\Eeq
\gianni{where we used \accorpa{primastimaA}{secondastimaA}, \pier{\eqref{hpdati} and \eqref{hpexA}} once more.}
Then, we multiply by $\dt\eta$ and integrate over~$Q_t$.
\pier{Thanks to the chain rule property (stated, e.g., in \cite[Lemme~3.3, p.~73]{Brezis})} \rev{and to the fact that $\eta(0) \in V$}, we obtain
\Beq
  \intQt |\dt\eta|^2
  + \frac \kappa 2 \iO |\nabla\eta(t)|^2
  + \rho \normaH{\eta(t)}
  = c (1+\rho)
  + \intQt g_2 \, \dt\eta
  \non
\Eeq
whence immediately
\Beq
  \norma{\dt\eta}_{\L2H}
  + \norma\eta_{\L\infty V}
  \leq c \bigl( \rho^{1/2} + 1 \bigr) .
  \label{terzastimaA}
\Eeq

\step
Fourth a priori estimate

\pier{%
We behave as we did for \eqref{secondastimaA}.
From \eqref{perterzaA} we have
\Beq
  - \kappa\Delta\eta(t) + \rho \sigma(t) = g_3(t) := g_2(t) - \dt\eta(t)
  \quad \aat .
  \non
\Eeq
Then, we formally multiply by $-\Delta\eta(t)$ and notice that
$\nabla\sigma(t)\cdot\nabla\eta(t)\geq0$ \aeO\
(at least formally; the inequality we need if $\Sign$ were replaced by $\Signeps$
would immediately follow from~\eqref{formulaSigneps}).
Hence, we get $\kappa^{1/2}\normaH{\Delta\eta(t)}\leq\normaH{g_3(t)}$  \aat.
By owing to~\eqref{perterzaA}, \eqref{terzastimaA} and elliptic regularity, we deduce~that
\Beq
  \norma\eta_{\L2W} 
  \leq c \bigl( \rho^{1/2} + 1 \bigr) .
  \label{quartastimaA}
\Eeq
}%

\step
\gianniold{Consequence}

\gianniold{Estimates \accorpa{primastimaA}{quartastimaA}
and assumption \eqref{hpexA}
imply for $\theta=\eta-\alpha\phi+\etastar$}
\Beq
  \gianniold{\norma\theta_{\L\infty H\cap \L2V} \leq c 
  \aand
  \norma\theta_{\H1H\cap\L2W}
  \leq c \bigl( \rho^{1/2} + 1 \bigr) \,.}
  \label{stimethetaA}
\Eeq

\step
Existence for Problem (A)

The above a~priori estimates are rigorous for the solution
to the approximating problem obtained by replacing
$\beta$ and $\Sign$ by the corresponding Yosida \regulariz ations.
Namely, one writes
\gianniold{%
\Beq
  \xi = \betaeps(\phi)
  \quad \aeQ
  \aand
  \sigma(t) = \Signeps(\eta(t))
  \quad \aat 
  \label{xisigmaeps}
\Eeq
}%
in place of \accorpa{terzaAbis}{quartaAbis}.
The approximating problem is more regular and has a solution 
$(\etaeps,\phieps,\xieps,\sigmaeps)$.
To see that, one can rewrite the approximating problem
by eliminating the time derivative $\dt\phi$ in \eqref{primaAbis}
on accout of~\eqref{secondaAbis}.
One obtains the Cauchy problem for a system of the form
\Beq
  \dt(\eta,\phi) + \calA(\eta,\phi) + \calB_\eps(\eta,\phi)
  = \calF
  \non
\Eeq
where $\calA$ is an unbounded operator in $\calH:=H\times H$,
$\calB_\eps:\calH\to\calH$ is a \Lip\ continuous perturbation
and $\calF$ is a source term.
Namely, $\calA$ acts as follows
\Bsist
  && \calA : (\eta,\phi) \mapsto \bigl( -\kappa\Delta\eta + \lambda \Delta\phi , - \nu\Delta\phi \bigr)
  \quad \hbox{for $(\eta,\phi) \in D(\calA) := W \times W$}  
  \non
  \\
  &&  \hbox{where} \quad
  \lambda := \kappa\alpha + (\ell-\alpha)\nu .
  \non
\Esist
Now, let us introduce the following inner product in~$\calH$
\Beq
  \bigl( (\eta,\phi) , (\tilde\eta,\tilde\phi) \bigr)_\calH
  := \iO \eta \gianniold{\tilde\eta} + \Bigl( \frac {\lambda^2}{\kappa\nu} + 1 \Bigr) \iO \gianniold\phi \tilde\phi .
  \non
\Eeq
Then, we have for $(\eta,\phi)\in D(\calA)$
\Bsist
  && \bigl( \calA(\eta,\phi) , (\eta,\phi) \bigr)_\calH
  = \iO \bigl( \kappa |\nabla\eta|^2 - \lambda \nabla\eta \cdot \nabla\phi + \frac {\lambda^2} \kappa \, |\nabla\phi|^2 \bigr)
  + \nu \iO |\nabla\phi|^2
  \non
  \\
  && \geq \frac \kappa 2 \iO |\nabla\eta|^2
  + \frac {\lambda^2} {2\kappa} \iO |\nabla\phi|^2
  + \nu \iO |\nabla\phi|^2
  \geq \frac \kappa 2 \iO |\nabla\eta|^2
  + \nu \iO |\nabla\phi|^2 .
  \non
\Esist
This shows that $\calA$ is monotone in $\calH$ with respect to that inner product.
Then, maximal monotonicy follows since the range of $\calA+\hbox{Id}_\calH$ 
is the whole of $\calH$ due to the Lax-Milgram theorem and elliptic regularity.
Therefore, the approximating problem has a solution
(see, e.g., \cite[Cor.~4.1 p.~181]{Show}).
So, by starting from the \analogue s of the above formal a priori estimates \rev{(that is, from the rigorous ones, for which properties \accorpa{dmoreau}{disG} have to be used)}
and \rev{owing to} standard weak, weakstar and strong compactness results
(see, e.g., \cite[Sect.~8, Cor.~4]{Simon}),
we have for a subsequence at least
\Bsist
  & \etaeps \to \eta
  & \hbox{weakly star in $\H1H\cap\L\infty V\cap\L2W$}
  \non
  \\
  && \hbox{and strongly in $\C0H$}
  \label{conveta}
  \\
  & \phieps \to \phi
  & \hbox{weakly star in $\H1H\cap\L\infty V\cap\L2W$}
  \non
  \\
  && \hbox{and strongly in $\C0H$}
  \label{convphi}
  \\
  & \gianniold{\xieps \to \xi}
  & \hbox{weakly in $\L2H$}
  \label{convchi}
  \\
  & \sigmaeps \to \sigma
  & \hbox{weakly star in $\L\infty H$} .
  \label{convsigma}
\Esist
\gianniold{We stress that $\xieps:=\betaeps(\phieps)$ and $\sigmaeps:=\Signeps(\etaeps)$, i.e.,
the same as in~\eqref{xisigmaeps}, where the subscripts $\eps$ were omitted for convenience.
Here, $\xi$~and $\sigma$ have the meaning given by \accorpa{convchi}{convsigma}.}
Clearly, the limits $\phi$, $\xi$ and $\sigma$ and the function $\theta$
computed from \eqref{defeta} satisfy the regularity requirements
and the estimates of the statement \gianniold{(see also~\eqref{stimethetaA})}.
Moreover, it follows that $\pi(\phieps)$ converges to $\pi(\phi)$ strongly in $\LQ2$
and that $\xi$ and $\sigma$ satisfy \accorpa{terzaAbis}{quartaAbis}
\pier{(because $\beta$ and $\Sign$ induce maximal monotone operators
on $\LQ2$ and $\L2H$, respectively,
and then they are weakly-strongly closed;
see, e.g., \cite[Cor.~2.4, p.~41]{Barbu})}.
Hence, $(\eta,\phi,\xi,\sigma)$ solves the original problem \PblAbis.

\step
Uniqueness for Problem (A)

We assume $\alpha=\ell$ and show that the solution is unique.
Let $(\eta_i,\phi_i,\xi_i,\sigma_i)$, $i=1,2$, be two solutions.
We write equations \accorpa{primaAbis}{secondaAbis} for both of them and take the differences.
If we set \gianniold{$\eta:=\eta_1-\eta_2$} and analogously define $\phi$, $\xi$ and~$\sigma$,
we obtain
\Bsist
  && \dt\eta - \kappa \Delta\eta + \kappa\ell \Delta\phi + \rho \sigma = 0
  \label{dprimaA}
  \\
  && \dt\phi - \nu \Delta\phi + \xi 
  = \gamma (\eta - \ell\phi) + \pi(\phi_2) - \pi(\phi_1) .
  \label{dsecondaA}
\Esist
Now, we multiply these equations by $\eta$ and $(\kappa\ell^2/\nu)\phi$, respectively,
sum up and integrate over~$Q_t$.
As $\pi$ is \Lip\ continuous, we~have
\Bsist
  && \frac 12 \iO |\eta(t)|^2
  + \frac {\kappa\ell^2} {2\nu} \iO |\phi(t)|^2
  + \kappa \intQt \bigl( |\nabla\eta|^2 - \ell \nabla\phi \cdot \nabla\eta + \ell^2 |\nabla\phi|^2 \bigr)
  \non
  \\
  && \quad {}
  + \rho \iot \bigl( \sigma(s) , \eta(s) \bigr)_H \, ds
  + \frac {\kappa\ell^2} \nu\intQt \xi \phi
  \leq c \intQt \bigl( |\eta|^2 + |\phi|^2 \bigr) .
  \non
\Esist
The last two terms on the \lhs\ are nonnegative by monotonicity
and the integral involving the gradients is estimated from below this way
\Beq
  \intQt \bigl( |\nabla\eta|^2 - \ell \, \nabla\phi \cdot \nabla\eta + \ell^2 |\nabla\phi|^2 \bigr)
  \geq \frac 12 \intQt \bigl( |\nabla\eta|^2 + \ell^2 |\nabla\phi|^2 \bigr).
  \label{frombelow}
\Eeq
At this point, we combine and apply the Gronwall lemma.
We conclude that $\eta=0$ and $\phi=0$, i.e.,
$\eta_1=\eta_2$ and $\phi_1=\phi_2$.
By comparison in \eqref{primaAbis} and \eqref{secondaAbis} written for both solutions,
we deduce that $\sigma_1=\sigma_2$ and \gianniold{$\xi_1=\xi_2$}, respectively.

\step
Further regularity

We assume \eqref{hppiuregA} and prove~\eqref{phipiureg}.
To this end, it suffices to perform the estimate corresponding to \eqref{phipiureg}
on the component $\phieps$ of the solution to the approximating problem sketched above, 
uniformly with respect to~$\eps$.
This can be done by a heavy calculation involving difference quotients.
Therefore, we confine ourselves to derive a formal a priori estimate.
We write equations \accorpa{secondaAbis}{terzaAbis}
by replacing $\beta$ by its Yosida \regulariz ation $\betaeps$ in the latter,
and formally differentiate with respect to time.
By writing $\phi$ instead of $\phieps$ for simplicity, we~have
\gianni{(\pier{see}~\eqref{primastimaA}, \eqref{terzastimaA}, and~\eqref{hpPi})}
\Beq
  \dt^2\phi - \nu \Delta\gianniold\dt\phi + \betaeps'(\phi) \dt\phi = g_3
  \quad \hbox{with} \quad
  \norma{g_3}_{\L2H} \leq c \bigl( \rho^{1/2} + 1 \bigr).
  \label{dtsecondaA}
\Eeq
Now, we multiply by~\gianniold{$\dt\phi$} and integrate over~$Q_t$.
We obtain
\Beq
  \frac 12 \iO |\dt\phi(t)|^2
  + \nu \intQt |\nabla\dt\phi|^2
  + \intQt \betaeps'(\phi) \, |\dt\phi|^2
  = \intQt g_3 \, \dt\phi
  + \frac 12 \iO |\dt\phi(0)|^2 .
  \non
\Eeq
As $\betaeps'$ is nonnegative by monotonicity, the only term that needs some \rev{treatment} is the last one on the \rhs.
We formally have from \eqref{secondaAbis}, the modified \eqref{terzaAbis} and the initial conditions
\Beq
  \dt\phi(0) 
  = \nu \Delta\phiz - \betaeps(\phiz) - \pi(\phiz) + \betti{\gamma}\thetaz \,.
  \label{dtphizA}
\Eeq
On the other hand, \eqref{defbetaz} implies that
$\normaH{\betaeps(\phiz)}\leq\normaH{\betaz(\phiz)}$.
Therefore, \gianniold{on account of~\eqref{hppiuregA},}
$\normaH{\dt\phi(0)}$ remains bounded and the estimate
\Beq
  \norma{\dt\phi}_{\L\infty H\cap\L2V}
  \leq c \bigl( \rho^{1/2} + 1 \bigr)
  \non
\Eeq
follows uniformly with respect to~$\eps$.
Thus, the same estimate holds for the limiting~$\phi$.
At this point, by comparison in~\eqref{secondaAbis},
we get a bound for the sum $-\nu\Delta\phi+\xi$ in $\L\infty H$
and the argument used to derive \eqref{secondastimaA}
\gianniold{(where $t$ is just a parameter)}
completes the regularity \eqref{phipiureg} and the estimate~\eqref{stimaApiureg}.
\gianniold{In order to conclude the proof of Theorem~\ref{WellposednessA},
we have to prove that the component $\theta$ of any solution
satisfying all the regularity requirements of the statement 
is bounded whenever we assume that $\thetaz\in\Linfty$ and $f\in\L\infty H$, in addition.
To this end, it suffices to write \eqref{primaA} in the form
\Beq
  \dt\theta - \kappa\Delta\theta
  = f - \rho\sigma - \ell\dt\phi
  \non
\Eeq
and observe that the \rhs\ of this equation belongs to~$\L\infty H$.
Then, we can argue, e.g., as in
\cite[Thm.~7.1, p.~181]{LSU} with $r=\infty$ and $q=2$,
where the case of the Dirichlet boundary conditions is treated in detail\rev{: by a careful check, the reader can make the necessary modifications to adapt the procedure to the case of the homogeneous} Neumann boundary conditions}.\QED

\subsection{\rev{Proof of Theorem~\ref{WellposednessB}}}

As the argument is quite similar to the previous one, we proceed quickly.
Also in this case, we introduce new unknowns and transform the problem.
\pier{Let us recall the} assumption \eqref{hpexB} on $\phistar$ and~set
\Beq
  \eta := \theta + \ell\phi , \quad
  \chi := \phi - \phistar
  \aand
  \xistar := \betaz(\phistar) .
  \label{defetachi}
\Eeq
Then, $\eta$ and $\chi$ have to satisfy the \analogue\ of \eqref{regsol}
and the new problem is the following
\Bsist
  && \dt\eta - \kappa \Delta\eta + \kappa\ell \Delta\chi
  = \gianniold{f - \kappa\ell \Delta\phistar}
  \quad \aeQ
  \label{primaBbis}
  \\
  && \dt\chi - \nu \Delta\chi + \xi - \xistar + \pi(\chi+\phistar)
  \non
  \\
  && = \gianniold{\gamma (\eta-\ell\chi-\ell\phistar) + \nu \Delta\phistar - \xistar - \rho \sigma}
  \quad \aeQ
  \qquad
  \label{secondaBbis}
  \\
  && \xi \in \beta(\chi+\phistar) \quad \aeQ
  \label{terzaBbis}
  \\
  && \sigma(t) \in \Sign(\chi(t))
  \quad \aat 
  \label{quartaBbis}
  \\
  && \eta(0) = \thetaz + \ell\phiz 
  \aand
  \chi(0) = \phiz - \phistar .
  \label{cauchyBbis}
\Esist

\step
First a priori estimate

We multiply \eqref{primaBbis} by $\eta$
and \eqref{secondaBbis} by $(\kappa\ell^2/\nu)\chi$,
integrate over $Q_t$ and sum up.
Then, we rearrange a little and use the \Lip\ continuity of $\pi$ and the Young inequality. 
\pier{Using also \eqref{hpexB}, we obtain}
\Bsist
  && \frac 12 \iO |\eta(t)|^2
  + \frac {\kappa\ell^2} {2\nu} \iO |\chi(t)|^2
  + \kappa \intQt
    \Bigl(
      |\nabla\eta|^2 - \ell \, \nabla\eta \cdot \nabla\chi + \ell^2 |\nabla\chi|^2
    \Bigr)
  \non
  \\
  && \quad {}
  + \frac {\kappa\ell^2} \nu\intQt (\xi-\xistar) \chi
  + \frac {\kappa\ell^2\gianniold\rho} \nu \iot \bigl( \sigma(s) , \chi(s) \bigr)_H \, ds
  \non
  \\
  && \leq c \intQt \bigl( |\eta|^2 + |\chi|^2 + 1 \bigr) .
  \non
\Esist
Now, we observe that \eqref{frombelow} \pier{can be applied}
and that the last two terms on the above \lhs\ are nonnegative by monotonicity.
Thus, by applying the Gronwall lemma,
we conclude that
\Beq
  \norma\eta_{\L\infty H\cap\L2V}
  + \norma\chi_{\L\infty H\cap\L2V}
  \leq c \,.
  \label{primastimaB}
\Eeq

\step
Second a priori estimate

We observe that \eqref{secondaBbis} looks like
\Beq
  \dt\chi - \nu\Delta\chi + \xi + \rho\sigma = \gianniold{g_1}
  \quad \hbox{with} \quad
  \norma{\gianniold{g_1}}_{\L2H} \leq c \,.
  \non
\Eeq
Therefore, multiplication by $\dt\chi$ and integration over $Q_t$ yield
\Bsist
  && \intQt |\dt\chi|^2
  + \frac {\betti{\nu}}{2} \iO |\nabla\chi(t)|^2
  + \iO \Beta(\chi(t)+\phistar)
  + \rho \normaH{\chi(t)}
  \non
  \\
  && \leq c (1+\rho)
  + \intQt \gianniold{g_1} \dt\chi
  \leq c \, (1+\rho)
  + \frac 12 \intQt |\dt\chi|^2 .
  \non
\Esist
We immediately deduce that
\Beq
  \norma\chi_{\H1H\cap\L\infty V} \leq c \bigl( \rho^{1/2} + 1 \bigr) .
  \label{secondastimaB}
\Eeq

\step
Further a priori estimates

\gianniold{We want to obtain}
\Beq
  \norma\eta_{\H1H\cap\L2W}
  + \norma\chi_{\H1H\cap\L2W}
  + \norma{\xi+\rho\sigma}_{\L2H}
  \leq c \bigl( \rho^{1/2} + 1 \bigr) .
  \label{stimeB}
\Eeq
To this end, we argue as we did for \accorpa{secondastimaA}{quartastimaA}
with just one modification of our argument concerning the pointwise estimate of~$\normaH{\Delta\chi(t)}$.
We still multiply by $\Delta\chi(t)$.
However, since $\phistar$ is not supposed to be a constant,
this requires some care and cannot be as simple as for~\eqref{secondastimaA}.
In order to be more precise on this delicate point, we consider the solution to the $\eps$-problem obtained by replacing
$\beta$ and $\Sign$ with their Yosida \regulariz ations $\betaeps$ and~$\Signeps$.
For simplicity, we avoid stressing the time $t$ for a while.
We write the \regulariz ed \eqref{secondaBbis} in the form
\Bsist
  && - \Delta\chi + \frac 1\nu \, \betaeps(\chi+\phistar) + \frac \rho\nu \, \Signeps\chi 
  = g_2
  \quad \hbox{with} 
  \non
  \\
  && g_2 := \frac 1\nu \bigl( - \dt\chi - \pi(\chi+\phistar) + \gamma (\eta-\ell\chi-\ell\phistar) + \nu \Delta\phistar \bigr)
  \non
\Esist
and read $-\Delta\chi$ as $-\Delta(\chi+\phistar)+\Delta\phistar$
when \pier{multiplying} the second term of the equation \pier{by~$-\Delta\chi$}.
\pier{Owing} to \eqref{formulaSigneps} (which also implies $\normaH{\Signeps\chi}\leq 1$)
and to the elementary inequalities~\eqref{elementary}, we obtain
\Bsist
  && \normaH{\Delta\chi}^2
  + \frac 1\nu \iO \betaeps'(\chi+\phistar) |\nabla(\chi+\phistar)|^2
  + \frac \rho\nu \iO \frac {|\nabla\chi|^2} {\max\graffe{\eps,\normaH\chi}}
  \non
  \\
  && = - \iO g_2 \Delta\chi
  - \frac 1\nu \iO \betaeps(\chi+\phistar) \Delta\phistar 
  \non
  \\
  && = - \iO g_2 \Delta\chi
  + \iO \Bigl( - \Delta\chi + \frac \rho\nu \Signeps\chi - g_2 \Bigr) \Delta\phistar 
  \non
  \\
  && = - \iO ( g_2 + \Delta\phistar ) \Delta\chi
  + \iO \Bigl( - g_2 + \frac \rho\nu \Signeps\chi \Bigr) \Delta\phistar 
  \non
  \\
  \separa
  && \leq (\normaH{g_2} + \normaH{\Delta\phistar}) \normaH{\Delta\chi}
  + \normaH{g_2} \, \normaH{\Delta\phistar}
  + \pier{\frac \rho\nu \, 
  \norma{\Delta\phistar}_{\rev H}}
  \non
  \\
  && \leq \frac 12 \, \normaH{\Delta\chi}^2
  + \frac 32 \, \normaH{g_2}^2
  + c \, (\rho + 1) \,.
  \non
\Esist
By recalling the meaning of~$g_2$, we conclude that we have \aat
\Bsist
  && \normaH{\Delta\chi(t)}^2
  \leq \frac 3 {\nu^2} \, \normaH{- \dt\chi(t) - \pi(\chi(t)+\phistar) + \gamma (\eta(t)-\ell\chi(t)-\ell\phistar) + \nu \Delta\phistar}^2
  \qquad
  \non
  \\
  && \phantom{\normaH{\Delta\chi(t)}^2\leq{}}
  + c \, (\rho + 1 ) \,.
  \label{stimaDeltachiB}
\Esist
Thus, the right bound for $\Delta\chi$ in $\L2H$ \gianni{follows from~\accorpa{primastimaB}{secondastimaB}}.
Then, $\xi+\rho\sigma$ is estimated in $\L2H$ by comparison in~\eqref{secondaBbis}
and the complete \eqref{stimeB} can be achieved like in the previous proof,
as said at the beginning.

\step
Existence for Problem (B)

One can proceed as for Problem (A).
Indeed, we have proved quite similar \gianniold{estimates}
\gianniold{(notice that \eqref{stimeB} also yields $\norma\xi_{\L2H}\leq c(\rho+1)$ 
since $\norma\sigma_{\L\infty H}\leq1$ by the definition of $\Sign$)
which} are completely rigorous when performed on the solution to the approximating problem
obtained by replacing $\beta$ and $\Sign$ by their Yosida \regulariz ations.
Moreover, the proof of the existence of a solution to the approximating problem is similar
to the one performed for Problem~(A).

\step
Uniqueness for Problem (B)

Let $(\eta_i,\phi_i,\xi_i,\sigma_i)$, $i=1,2$, be two solutions
and define $\eta_i$ and $\chi_i$ according to~\eqref{defetachi}.
By proceeding in the same way as we did for Problem~(A),
we easily obtain $\eta_1=\eta_2$ and $\chi_1=\chi_2$,
whence $\theta_1=\theta_2$ and $\phi_1=\phi_2$,
i.e., the first sentence of Theorem~\ref{WellposednessB} about uniqueness.
Now, assume $\beta$ to be single-valued.
Then, $\xi_1=\xi_2$ since $\phi_1=\phi_2$.
Finally, by comparison in~\eqref{secondaB},
we also deduce that $\sigma_1=\sigma_2$.
This concludes the proof of Theorem~\ref{WellposednessB}.\QED

\subsection{\rev{Proof of Theorem~\ref{WellposednessC}}}

As we did for Problem~(B),
we introduce the new unknowns $\eta$ and $\chi$ by means of~\eqref{defetachi}
and deal with the following new problem:
\Bsist
  && \dt\eta - \kappa \Delta\eta + \kappa\ell \Delta\chi
  = f \rev{{}- \kappa \ell \Delta \phistar}
  \quad \aeQ
  \label{primaCbis}
  \\
  && \dt\chi - \nu \Delta\chi + \xi - \xistar + \pi(\chi+\phistar)
  \non
  \\
  && = \gamma (\eta-\ell\chi-\ell\phistar) \gianniold{{}+\rev\nu\Delta\phistar} - \xistar - \rho \sigma
  \quad \aeQ
  \qquad
  \label{secondaCbis}
  \\
  && \xi \in \beta(\chi+\phistar) \quad \aeQ
  \label{terzaCbis}
  \\
  && \sigma \in \sign\chi
  \quad \aeQ
  \label{quartaCbis}
  \\
  && \eta(0) = \thetaz + \ell\phiz 
  \aand
  \chi(0) = \phiz - \phistar .
  \label{cauchyCbis}
\Esist

\step
Existence and uniqueness for Problem (C)

This problem differs from Problem~(B) just in~\eqref{quartaCbis},
where $\sign$ appears in place of the non-local operator~$\Sign$.
Therefore, both existence and partial uniqueness can be obtained by the same argument.

It remains to prove the regularity part of
Theorem~\ref{WellposednessC} and the estimates.
For the regularity of $\phi$ and~$\xi$, 
one could combine the techniques used for Problems~(A) and~(B), while the
regularity of $\theta$ is classical, on account of~\eqref{hppiuregC}
and of the regularity of $\dt\phi$ already proved.
However, the forthcoming argument also shows the desired regularity.
What needs much more care is the control of the constants entering \accorpa{chibddC}{thetabddC}.
\gianni{This forces us to perform a number of a~priori estimates
which we derive just formally, for brevity.}

\step
First a priori estimate

As we did for~\eqref{primastimaB}, we multiply \eqref{primaCbis} by $\eta$
and \eqref{secondaCbis} by $(\kappa\ell^2/\nu)\chi$,
integrate over $Q_t$ and sum up.
Then, we owe to \eqref{frombelow}, the \Lip\ continuity of~$\pi$,
the Young inequality and the Gronwall lemma.
We obtain
\Beq
  \norma\eta_{\L\infty H\cap\L2V}
  + \norma\chi_{\L\infty H\cap\L2V}
  \leq c \,.
  \label{primastimaC}
\Eeq

\step
Second a priori estimate

We write \eqref{secondaCbis} as
\Beq
  \dt\chi - \nu\Delta\chi + \xi + \rho\sigma = g_1
  \quad \hbox{with} \quad
  \norma{g_1}_{\L2H} \leq c 
  \non
\Eeq
and multiply \pier{it} by $\dt\chi$.
As for~\eqref{secondastimaB}, we have
\Beq
  \norma\chi_{\H1H\cap\L\infty V} \leq c \bigl( \rho^{1/2} + 1 \bigr) .
  \label{secondastimaC}
\Eeq

\step
Third a priori estimate

\gianniold{There holds 
\Bsist
  && \normaH{\Delta\chi(t)}^2
  \leq \frac 3 {\nu^2} \, \normaH{- \dt\chi(t) - \pi(\chi(t)+\phistar) + \gamma (\eta(t)-\ell\chi(t)-\ell\phistar) + \nu \Delta\phistar}^2
  \qquad
  \non
  \\
  && \phantom{\normaH{\Delta\chi(t)}^2\leq{}}
  + c \, (\rho + 1 ) 
  \label{stimaDeltachiC}
\Esist
i.e., the same as \eqref{stimaDeltachiB}.
Inequality \eqref{stimaDeltachiC} can be proved with the same calculations that led to~\eqref{stimaDeltachiB}
with obvious changes in the proof 
(like $\normaH{\signeps\chi}\leq|\Omega|^{1/2}$ in place of $\normaH{\Signeps\chi}\leq1$,
whence just a different value of the final~$c$).
From this and the previous estimates, we deduce that}
\Beq
  \norma{\Delta\chi}_{\L2H}
  \leq c \bigl( \rho^{1/2} + 1 \bigr).
  \label{terzastimaC}
\Eeq

\step
Fourth a priori estimate

Now, we multiply \eqref{primaCbis} by~$\dt\eta$, integrate over~$Q_t$ and get
\Beq
   \intQt |\dt\eta|^2
  + \frac \kappa 2 \iO |\nabla\eta(t)|^2
  = \gianni{\frac \kappa 2 \iO |\nabla\eta(0)|^2
  + \intQt (f - \kappa\ell \Delta\chi \rev{{}- \kappa \ell \Delta \phistar}) \dt\eta}.
  \non
\Eeq
\pier{From \eqref{hpdati} and \eqref{stimaDeltachiC}, we} immediately infer that
\Beq
  \norma\eta_{\H1H\cap\L\infty V}
  \leq c \bigl( \rho^{1/2} + 1 \bigr).
  \label{quartastimaC}
\Eeq

\step
Fifth a priori estimate

\gianni{We start from~\eqref{secondaCbis}
and smooth the monotone nonlinearities
by replacing them with their Yosida approximations.
\pier{By differentiating with respect to time}, we~have}
\Beq
  \dt^2\chi - \nu \Delta\dt\chi + \bigl\{ \betaeps'(\chi+\phistar) + \rho \signeps'(\chi) \bigr\} \dt\gianniold\chi
   = g_3
  \label{dtsecondaC}
\Eeq
\pier{where
$g_3 := \gamma (\dt\eta - \ell\dt\chi) - \pi'(\chi+\phistar) \dt\chi$,
and we can read the initial value}
\Beq
  \dt\chi(0) 
  = \nu \Delta\phiz - \betaeps(\phiz) - \pi(\phiz) + \gamma\thetaz - \rho \signeps(\phiz-\phistar) .
  \label{dtphizC}
\Eeq
Notice that
\Beq
  \norma{g_3}_{\L2H} \leq c \bigl( \rho^{1/2} + 1 \bigr)
  \aand
  \normaH{\dt\chi(0)}
  \leq \rho \, |\Omega|^{1/2} + c 
  \label{stimerhsC}
\Eeq
\pier{thanks to \eqref{hppiuregA}}.
Thus, by multiplying \eqref{dtsecondaC} by~$\dt\gianni\chi$, integrating over~$Q_t$,
observing that $\betaeps'$ and $\signeps'$ are nonnegative,
and using~\eqref{secondastimaC} for $\dt\chi$ and~\eqref{stimerhsC},
we obtain
\Beq
  \frac 12 \iO |\dt\chi(t)|^2
  + \nu \intQt |\nabla\dt\chi|^2
  \leq \intQt |g_3| \, |\dt\chi|
  + \frac 12 \iO |\dt\chi(0)|^2 
  \leq \frac {\rho^2 |\Omega|} 2 
  + c \, \rho + c \,.
  \non
\Eeq
We deduce the following estimates for~$\dt\chi$
\Beq
  \norma{\dt\chi}_{\L\infty H}^2
  \leq \rho^2 |\Omega| + c \, \rho + c 
  \aand
  \norma{\nabla\dt\chi}_{\L2H}^2
  \leq \frac {\rho^2 |\Omega|} {2\nu} + c \, \rho + c \,.
  \label{quintastimaC}
\Eeq

\step
Sixth a priori estimate 

\gianniold{We use~\eqref{primastimaC}, \eqref{stimaDeltachiC} and \eqref{quintastimaC}.}
We deduce that
\Beq
  \norma{\Delta\chi}_{\L\infty H}^2
  \leq \frac {6\rho^2 \, |\Omega|} {\gianni{\nu^2}} + c \, \rho + c \,.
  \non
\Eeq
Now, \gianniold{we recall the definition \eqref{defnormaW} of~$\normaW\cpto$}.
Thus, we also have
\Beq
  \norma\chi_{\L\infty W}^2
  \leq \gianniold{\frac {6\rho^2 \, |\Omega|^{7/3}} {\gianni{\nu^2}}} + c \, \rho + c \,.
  \non
\Eeq
Finally, we apply~\eqref{embedding}.
We conclude that $\chi\in\LQ\infty$ and that
\Beq
  \norma\chi_\infty^2
  \leq \gianniold{\frac {6\rho^2 \, \CO^2 \,  |\Omega|^{7/3}} {\gianni{\nu^2}}} + c \, \rho + c 
  \non
\Eeq
whence also \gianniold{(by the first elementary inequality \eqref{elementary})}
\Beq
  \norma\chi_\infty
  \leq \frac {6^{1/2} \, \rho \, \CO \,  |\Omega|^{7/6}} {\gianni\nu} + c \bigl( \rho^{1/2} + 1 \bigr) 
  \leq 2 \, \frac {6^{1/2} \, \rho \, \CO \,  |\Omega|^{7/6}} {{\gianni\nu}} + c  \,.
  \label{sestastimaC}
\Eeq
Hence, if we \gianniold{choose} the last value of~$c$ as~$C_6$,
we see that \eqref{chibddC} holds with $\CS$ as in~\eqref{costC}.

\step
Seventh a priori estimate

On account of the \pier{regularity of $f$ in}~\eqref{hppiuregC},
we formally differentiate \eqref{primaCbis} with respect to time 
and test the resulting equation by~$\dt\eta$.
As $\dt\eta(0)$, which \pier{is recovered} from~\eqref{primaCbis},
is bounded in~$H$ by a constant due to \eqref{hppiuregC},
we obtain \gianni{by~\eqref{quartastimaC}}
\Bsist
  && \frac 12 \iO |\dt\eta(t)|^2
  + \kappa \intQt |\nabla\dt\eta|^2
  = \frac 12 \iO |\dt\eta(0)|^2
  + \kappa\ell \intQt \nabla\dt\eta \cdot \nabla\dt\chi
  \gianni{{}+ \intQt \dt f \, \dt\eta}
  \non
  \\
  && \leq c + \kappa \intQt |\nabla\dt\eta|^2
  + \frac {\kappa\ell^2} 4 \intQ |\nabla\dt\chi|^2 .
  \non
\Esist
\pier{Owing} to the second of~\eqref{quintastimaC}, we infer that
\Beq
  \norma{\dt\eta}_{\L\infty H}^2
  \leq \frac {\kappa\ell^2\rho^2|\Omega|} {4\nu} + c \, \rho + c \,.
  \label{settimastimaC}
\Eeq

\step
Eighth a priori estimate

\gianni{By recalling that $\theta=\eta-\ell\chi-\ell\phistar$ by~\eqref{defetachi},
the first \pier{inequality in} \eqref{quintastimaC} and \pier{estimate} \eqref{settimastimaC} yield}
\Beq
  \norma{\dt\theta}_{\L\infty H}
  \leq \hat C \, \rho \, |\Omega|^{1/2}
  + c \rho^{1/2} + c 
  \quad \hbox{where} \quad
  \hat C := \frac {\kappa^{1/2}\ell} {2\nu^{1/2}} + \ell \,.
  \label{stimadtthetaC}
\Eeq
Once such an estimate is obtained,
we can recover a bound for $\Delta\theta$ from \eqref{primaC}
and repeat for $\theta$ what we have done for~$\rev{\chi}$.
Here is the quick sequence of deductions.
First, we have
\Bsist
  && \norma{\Delta\theta}_{\L\infty H}
  \leq \frac 1\kappa
  \bigl(
    \norma f_{\L\infty H}
    + \norma{\dt\theta}_{\L\infty H}
    + \ell \norma{\dt\chi}_{\L\infty H}
  \bigr)
  \non
  \\
  && \leq \frac {\hat C+\ell} \kappa \, \rho |\Omega|^{1/2}
  + c \rho^{1/2} + c 
  \non
\Esist
and we derive 
\Bsist
  && \norma\theta_{\L\infty W}^2
  \leq \norma\theta_{\L\infty H}^2
  + |\Omega|^{4/3} \norma{\Delta\theta}_{\L\infty H}^2
  \non
  \\
  && \leq 4 \, \frac {(\hat C+\ell)^2} {\kappa^2} \, \rho^2 |\Omega|^{7/3}
  + c (\rho + 1) \,.
  \non
\Esist
\gianni{Therefore}
\Beq
  \norma\theta_\infty
  \leq \CO \, 2 \, \frac {\hat C+\ell} \kappa \, \rho |\Omega|^{7/6}
  + c \bigl( \rho^{1/2} + 1 \bigr)
  \leq \CO \, 4 \, \frac {\hat C+\ell} \kappa \, \rho |\Omega|^{7/6}
  + c 
  \non
\Eeq
so that \eqref{thetabddC} holds with the last value of~$c$ as~$C_7$
and $\CS$ as in~\eqref{costC}
\gianni{(recall the value of $\hatC$ in~\eqref{stimadtthetaC})}.


\section{Existence of sliding modes}
\label{SLIDINGMODES}
\setcounter{equation}{0}

This section is devoted to the proof of \gianniold{Theorems~\ref{SlidingA}, \ref{SlidingB} and~\ref{SlidingC}}.
The argument we use to prove the existence of sliding modes \gianniold{in the first two cases}
relies on the following lemma, 
which ensures the existence of an extinction time $\Tstar$ for a real function.

\Blem
\label{Extinction}
Let $\az,\bz,\psiz,\rho\in\erre$ be such that
\Beq
  \az, \, \bz, \, \psiz \geq 0
  \aand
  \rho > \az^2 + 2\bz + 2 \, \frac \psiz T
  \label{hpext}
\Eeq
and let $\psi:[0,T]\to[0,+\infty)$ be an absolutely continuous function
satisfying $\psi(0)=\psiz$ and
\Beq
  \psi' + \rho \leq \az \, \gianniold{\rho^{1/2}} + \bz
  \quad \hbox{a.e.\ in the set \ $P:=\{t\in(0,T):\ \psi(t)>0\}$}.
  \label{diseqa}
\Eeq
Then, the following conclusions hold true.

\noindent
$i)$~If $\psiz=0$, then $\psi$ vanishes identically.

\noindent
$ii)$~If $\psiz>0$, there exists $\Tstar\in(0,T)$ satisfying
$\Tstar\leq2\psiz/(\rho-\az^2-2\bz)$ such that
$\psi$ is strictly decreasing in $(0,\Tstar)$ and $\psi$ vanishes in~$[\Tstar,T]$.
\Elem

\Bdim
Assumption \eqref{diseqa} and the Young inequality imply that
\Beq
  \psi' \leq - s_0
  \quad \hbox{a.e.\ in $P$},
  \quad \hbox{where} \quad
  s_0 := \frac 12 \, \rho - \frac 12 \, \az^2 - \bz 
  \label{diseqb}
\Eeq
and we notice that \eqref{hpext} implies
\Beq
  s_0 > \frac \psiz T \,.
  \label{slopeok}
\Eeq
In particular, $s_0>0$.
Moreover,
if $0\leq t_1<t_2\leq T$ and $(t_1,t_2)\subseteq P$, then
\Beq
  \psi(t_1) 
  = \psi(t_2) - \int_{t_1}^{t_2} \gianniold{\psi'}(t) \, dt
  \geq \psi(t_2) + s_0 (t_2-t_1)
  \geq s_0 (t_2-t_1) 
  > 0 .
  \label{psipos}
\Eeq
Now, we prove the lemma.

$i)$~By contradiction, let $P$ be non-empty.
So, we can pick a connected component of~it.
This is an open interval~$(a,b)$
and we can apply \eqref{psipos} to obtain $\psi(a)>0$.
Thus, $a>0$ since $\psiz=0$, whence
$\psi>0$ also in $(a',a]$ for some $a'<a$.
This contradicts the definition of connected component.

$ii)$~As $\psiz>0$, we can define the strictly positive number $\Tstar$ by setting
\Beq
  \Tstar := \sup \{ t\in(0,T): \ \psi(s)>0 \ \hbox{for every $s\in(0,t)$} \} .
  \non
\Eeq
By \eqref{psipos} with $t_1=0$ and $t_2=\Tstar$ and~\eqref{slopeok}, we have $\psiz\geq s_0\Tstar$, 
whence $\Tstar\leq\psiz/s_0<T$, i.e.,
the first conditions of the statement.
Furthermore, $\psi'\leq-s_0<0$ in $(0,\Tstar)$, 
so that $\psi$ is strictly decreasing in this interval.
Finally, we have to show that $\psi$ vanishes in $[\Tstar,T]$
and we argue by contradiction by assuming that $P\cap(\Tstar,T)\not=\emptyset$
and picking a connected component of this set.
This is an open interval $(a,b)$, with $\Tstar\leq a<b\leq T$, in principle.
However, $a=\Tstar$ would contradict the definition of~$\Tstar$, whence $a>\Tstar$.
Therefore, by applying \eqref{psipos} with $t_1=a$ and $t_2=b$, 
we obtain $\psi(a)>0$ and the definition of connected component is contradicted
as in the previous case.
\Edim

\step
Proof of Theorem~\ref{SlidingA}

Let $(\theta,\phi,\xi,\sigma)$ be a solution to problem \PblA\
\gianni{given by \accorpa{conveta}{convsigma}}.
We show that this solution \rev{fulfills} the requirements of the statement.
First of all, we observe that estimates \accorpa{stimaAuno}{stimaAdue} and \eqref{stimaApiureg}
hold for the approximating solution, by construction.
Moreover, $f\in\L\infty H$ by assumption.
Hence, we can write the modified \eqref{primaAbis} in the form
\Bsist
  && \dt\etaeps - \kappa \Delta\etaeps + \rho \sigmaeps
  = \geps 
  := f - (\ell-\alpha) \dt\phieps - \kappa\alpha \Delta\gianniold\phieps \pier{{}+{}} \kappa \Delta\etastar
  \label{primaAter}
  \\
  && \norma\geps_{\L\infty H} \leq C \bigl( \rho^{1/2} + 1 \bigr) 
  \label{rhsAter}
\Esist
where $C$ depends only on the structure and the data involved in the statement.
At this point, we~set
\Beq
  \rhostar := C^2 + 2C + \frac 2T \, \normaH{\thetaz+\alpha\phiz-\etastar}
  \label{defrhostarA}
\Eeq
and assume $\rho>\rhostar$.
We also set
\Beq
  \psi(t) :=\normaH{\eta(t)}
  \aand
  \psieps(t) :=\normaH{\etaeps(t)}
  \quad \hbox{for $t\in[0,T]$}.
  \label{defpsiA}
\Eeq
\gianniold{Now}, by assuming $h\in(0,T)$ and $t\in(0,T-h)$,
we multiply \eqref{primaAter} by $\sigmaeps$ and integrate over
$(t,t+h)\times\Omega$.
We obtain
\Bsist
  && \ioth \bigl( \dt\etaeps(s) , \sigmaeps(s) \bigr)_H \, ds
  + \kappa \ioth \iO \nabla\etaeps \cdot \nabla\sigmaeps
  + \rho \ioth \normaH{\sigmaeps(s)}^2 \, ds
  \non
  \\
  && = \ioth \bigl( \geps(s) , \sigmaeps(s) \bigr)_H \, ds .
  \label{slidA}
\Esist
As \eqref{moreau} and \eqref{dmoreau} imply that
\Beq
  \bigl( \dt\etaeps(t) , \sigmaeps(t) \bigr)_H 
  = \frac d {dt} \int_0^{\psieps(t)} \min \{ s/\eps , 1 \} \, ds
  \quad \aat
  \non
\Eeq
we have for the first term of \eqref{slidA}
\Beq
  \ioth \bigl( \dt\etaeps(s) , \sigmaeps(s) \bigr)_H \, ds
  = \int_{\psieps(t)}^{\psieps(t+h)} \min \{ s/\eps , 1 \} \, ds .
  \non
\Eeq
The second integral in \eqref{slidA} is nonnegative.
Indeed, \eqref{formulaSigneps} implies
\Beq
  \nabla\etaeps(t) \cdot \nabla\sigmaeps(t)
  = \frac {|\nabla\etaeps(t)|^2} {\max\{\eps,\normaH{\etaeps(t)}\pier\}} 
  \geq 0
  \quad \aeO \pier, \enskip \aat .
  \non
\Eeq
As $\normaH{\sigmaeps(s)}\leq 1$ for every $s$ and \eqref{rhsAter} holds, 
we deduce from~\eqref{slidA}
\Beq
  \int_{\psieps(t)}^{\psieps(t+h)} \min \{ s/\eps , 1 \} \, ds
  + \rho \ioth \normaH{\sigmaeps(s)}^2 \, ds
  \leq h \, C \bigl( \rho^{1/2} + 1 \bigr) .
  \non
\Eeq
At this point, we let $\eps$ tend to zero.
As we are assuming that \eqref{conveta} and \eqref{convsigma} hold at least for a subsequence,
we infer~that
\Bsist
  && \psi(t+h) - \psi(t)
  + \rho \ioth \normaH{\sigma(s)}^2 \, ds
  \non
  \\
  && \leq \lim_{\eps\seto 0} \int_{\psieps(t)}^{\psieps(t+h)} \min \{ s/\eps , 1 \} \, ds
  + \rho \liminf_{\eps\seto 0} \ioth \normaH{\sigmaeps(s)}^2 \, ds
  \leq h \, C \bigl( \rho^{1/2} + 1 \bigr) 
  \non
\Esist
for every $h\in(0,T)$ and $t\in(0,T-h)$.
This implies that
\Beq
  \psi'(t)
  + \rho \normaH{\sigma(t)}^2
  \leq C \bigl( \rho^{1/2} + 1 \bigr)
  \quad \aat .
  \non
\Eeq
As $\normaH{\sigma(t)}=1$ if $\norma{\eta(t)}>0$ by~\eqref{signv},
we can apply the lemma with $\az=\bz=C$
and we observe that our condition $\rho>\rhostar$ completely fits the assumptions
by~\eqref{defrhostarA}.
Thus, we find $\Tstar\in\gianniold[0,T)$ such that $\eta(t)=0$ for every $t\in[\Tstar,T]$,
i.e., \eqref{slidingA}.\QED

\step
Proof of Theorem~\ref{SlidingB}

By arguing as in the previous proof,
we pick a solution $(\theta,\phi,\xi,\sigma)$ to problem \PblB\
obtained as the limit of the solution $(\thetaeps,\phieps,\xieps,\sigmaeps)$
of the corresponding approximating problem
and show that all the requirements of the statement are fulfilled.
We introduce the functions $\eta$ and $\chi$ defined by \eqref{defetachi}
and the \analogue s $\etaeps$ and $\chieps$, 
and owe to \accorpa{stimaBuno}{stimaBdue}
for the approximating solution.
Therefore, we can rewrite \pier{the equation approximating} \eqref{secondaBbis} in the form
\Bsist
  && \dt\chieps - \nu \Delta\gianniold\chieps + \betaeps(\chieps+\phistar) - \betaeps(\phistar) + \rho \gianniold\sigmaeps
  \non
  \\
  && = \geps := \gianniold{\gamma (\etaeps-\ell\chieps-\ell\phistar) \gianniold{{}+\nu\Delta\phistar} - \betaeps(\phistar) \rev{{}- \pi (\chieps +\phistar)}}\qquad
  \label{secondaBter}
\Esist
\rev{with}
\Beq
   \norma\geps_{\L\infty H} \leq C  
  \label{rhsBter}
\Eeq
where $C$ depends only on the structure and the data involved in the statement.
At this point, we~set
\Beq
  \rhostar := 2C + \frac 2T \, \normaH{\phiz-\phistar}
  \label{defrhostarB}
\Eeq
and assume $\rho>\rhostar$.
We also set
\Beq
  \psi(t) :=\normaH{\chi(t)}
  \aand
  \psieps(t) :=\normaH{\chieps(t)}
  \quad \hbox{for $t\in[0,T]$}.
  \label{defpsiB}
\Eeq
\gianniold{Now, we} multiply \eqref{secondaBter} by $\sigmaeps$ and integrate over~$(t,t+h)\times\Omega$ as before.
We obtain 
\Bsist
  && \ioth \bigl( \dt\chieps(s) , \sigmaeps(s) \bigr)_H \, ds
  + \nu \ioth \iO \nabla\chieps \cdot \nabla\sigmaeps
  \non
  \\
  && \quad {}
  + \ioth \bigl( \betaeps(\chieps(s)+\phistar) - \betaeps(\phistar) , \Signeps(\chieps(s)) \bigr)_{\!\rev H} \, ds
  + \rho \ioth \normaH{\sigmaeps(s)}^2 \, ds
  \non
  \\
  && = \ioth \bigl( \geps(s) , \sigmaeps(s) 
  \bigr)_{\!\rev H} \, \gianniold{ds}.
  \label{slidB}
\Esist
The first two terms and the \gianni{\lhs}\ can be dealt with as in the previous proof.
The third integral on the \lhs\ is nonnegative since
the two factors of the product have the same sign.
Therefore, by arguing as above and then \gianniold{letting} $\eps$ tend to zero, we~obtain
\Beq
  \psi'(t)
  + \rho \normaH{\sigma(t)}^2
  \leq C .
  \non
\Eeq
As $\normaH{\sigma(t)}=1$ if $\norma{\chi(t)}>0$ by~\eqref{signv},
we can apply the lemma with $\az=0$ and $\bz=C$
since $\rho>\rhostar$ (see~\eqref{defrhostarB}).
Thus, we find $\Tstar\in\gianniold[0,T)$ such that $\chi(t)=0$ for every $t\in[\Tstar,T]$.
This condition coincides with~\eqref{slidingB}.\QED

\step
Proof of Theorem~\ref{SlidingC}

\gianniold{%
For Problem~(C) we use a different argument
since we cannot apply Lemma~\ref{Extinction}.
Our method relies on a comparison \rev{technique}
on $\chi:=\phi-\phistar$,
where $(\theta,\phi,\xi,\sigma)$ is the solution we are dealing with,
\gianni{by introducing the solution $w$ of an ordinary Cauchy problem with a well-chosen \rhs}.
\gianni{The function $\chi$} has the same regularity of $\phi$ and satisfies
\begin{align}
  & \dt\chi - \nu \Delta\chi + \xi - \xistar + \pi(\chi+\phistar) + \rho \, \sigma
  = \gamma \theta + \nu\Delta\phistar - \xistar
  \quad \pier{\aeQ}
  \label{secondachi}
  \\
  & \hbox{where} \quad
  \xistar := \betaz(\phistar)
  \label{defxistar}
  \\
  & \xi \in \beta(\chi+\phistar) 
  \aand
  \sigma \in \sign \chi
  \quad \pier{\aeQ}
  \label{terzachi}
  \\
  & \dn\chi = 0
  \quad \pier{\aeS}
  \aand 
  \chi(0) = \chiz := \phiz - \phistar .
  \label{bicchi}
\end{align}
Our starting point is just estimate~\eqref{thetabddC}, 
i.e., we only \gianni{suppose} that the constants $\CS$, $\CO$ and $C_7$ satisfy~it
and do not require that they are constructed as in the proof of Theorem~\ref{WellposednessC}.
\gianni{In order to introduce the ingredients of the Cauchy problem mentioned above,
we set for convenience}
\Bsist
  && \Mz := \norma\chiz_\infty 
  \aand
  \Mpi := L(\Mz+\norma\phistar_\infty) + |\pi(0)|
  \label{defMz}
  \\
  && \Mrho := \rho \, \CS\CO |\Omega|^{7/6} + C_7 
  \label{defMrho}
  \\
  && \Arho := \gamma \Mrho + \nu \norma{\Delta\phistar}_\infty + \norma\xistar_\infty + \Mpi
  \label{defArho}
\Esist
\gianni{where $L$ is the \Lip\ constant of~$\pi$}.
\gianni{We observe that (cf.~\eqref{thetabddC})}
\Beq
  \norma\theta_\infty \leq \Mrho
  \aand
  |\pi(\phistar\pm r)| \leq \Mpi
  \quad \aeO
  \quad \hbox{for every $r\in [0,\Mz] $}.
  \label{thetapibdd}
\Eeq
We assume \eqref{smallness}
and define $\rhostar$ as the solution to $\rho=\Arho+\Mz/T$, i.e.,
\Beq
  \rhostar := \frac
  {\gamma C_7 + \nu \norma{\Delta\phistar}_\infty + \norma\xistar_\infty + \Mpi + \Mz/T}
  {1 - \gamma \, \CS\CO |\Omega|^{7/6}} \,.
  \label{defrhostarC}
\Eeq
We claim that $\rhostar$ \rev{fulfills} the properties of the statement.
So, we fix $\rho>\rhostar$ and consider any solution 
\pier{of the transformed problem according}
to \accorpa{secondachi}{bicchi}.
We observe that our assumption $\rho>\rhostar$ implies
\Beq
  \rho > \Arho + \frac \Mz T \,,
  \quad \hbox{whence also} \quad
  \rho > \Arho 
  \label{perTstar}
\Eeq
and \pier{we can} set
\Beq
  \Tstar := \frac \Mz {\rho-\Arho} \,.
  \label{defTstarC}
\Eeq
\pier{Hence,} \eqref{perTstar} ensures that the definition of $\Tstar$ is meaningful
and that $\Tstar\geq0$.
More precisely, $\Tstar=0$ if $\Mz=0$, i.e., $\phiz=\phistar$,
and $\Tstar>0$ otherwise.
The first \pier{inequality in} \eqref{perTstar} implies that $\Tstar<T$.
The rest of the proof is devoted to prove that $\chi(t)=0$ for every $t\in[\Tstar,T]$.
This is done by comparison arguments, as \gianni{mentioned} at the beginning.
We introduce the ordinary Cauchy problem
\Beq
  w' + \rho \zeta = \Arho, \quad
  \zeta \in \sign w
  \aand
  w(0) = \Mz \,.
  \label{odew}
\Eeq
As $\Arho/\rho\in[0,1)\subset\sign0$ by \eqref{perTstar}, 
one checks that its unique solution is given by
\Beq
  w(t) = \bigl( \Mz - (\rho-\Arho) t \bigr)^+
  \quad \hbox{for $t\in[0,T]$}.
  \label{soluzw}
\Eeq
Notice that $0\leq w\leq\Mz$ and that $w$ vanishes on $[\Tstar,T]$ by the definition \eqref{defTstarC} of~$\Tstar$.
Thus, by also reading $w$ as a space independent function defined in $Q$ rather than in~$(0,T)$,
it suffices to prove that $|\chi|\leq w$ \aeQ.
To this end, we observe that $w$ trivially satisfies the homogeneous Neumann boundary condition
and write \eqref{odew} in the following forms
\Bsist
  && \dt w - \nu \Delta w + \pi(\phistar+w) + \rho\zeta
  = \Arho + \pi(\phistar+w)
  \label{odewa}
  \\
  && \dt w - \nu \Delta w - \pi(\phistar-w) - \rho (-\zeta)
  = \Arho - \pi(\phistar-w)
  \label{odewb}
  \\
  && \hbox{with} \quad
  \zeta \in \sign w
  \quad \hbox{or, equivalently,} \quad
  -\zeta \in \sign(-w) \,.
  \non
\Esist
We set $\psi:=(\chi-w)^+$, the positive part of~$\chi-w$, and multiply the difference 
between \eqref{secondachi} and \eqref{odewa} by~$\psi$.
By accounting for \eqref{thetapibdd} and the definition \eqref{defArho} of~$\Arho$, we have
\Bsist
  && \frac 12 \iO |\psi(t)|^2
  + \nu \intQt |\nabla\psi|^2
  + \intQt (\xi-\xistar) \psi
  + \rho \intQt (\sigma-\zeta) \psi
  \non
  \\
  && \quad {}
  + \intQt \bigl( \pi(\phistar+\chi) - \pi(\phistar+w) \bigr) \psi
  \non
  \\\
  && = \intQt \bigl( \gamma\theta + \nu\Delta\phistar - \xistar - \Arho - \pi(\phistar+w) \bigr) \psi
  \non
  \\
  && \leq \intQt \bigl( \gamma\Mrho + \nu \norma{\Delta\phistar}_\infty + \norma\xistar_\infty - \Arho + \Mpi \bigr) \psi
  = 0 .
  \non
\Esist
Now, we observe that the integrals on the \lhs\ involving $\xi$ and $\sigma$ are nonnegative:
indeed, where $\psi\not=0$, we have $\psi>0$ and $\chi>w$,
whence $\phistar+\chi>\phistar+w\geq\phistar$, so that $\xi\geq\xistar$ and $\sigma\geq\zeta$.
On the other hand, we~have
\Beq
  \intQt \bigl( \pi(\phistar+\chi) - \pi(\phistar+w) \bigr) \psi
  \geq - L \intQt |\chi-w| \psi
  = - L \intQt |\psi|^2 .
  \non
\Eeq
Therefore, we deduce that
\Beq
  \iO |\psi(t)|^2 \leq L \intQt |\psi|^2 .
  \label{lastgronwall}
\Eeq
By applying the Gronwall lemma, we conclude that $\psi=0$, i.e., $\chi\leq w$.
Now, we set $\psi:=(\chi+w)^-$, the negative part of~$\chi+w$, add equations \eqref{secondachi} and \eqref{odewb} to each other
and multiply the resulting equality by~$-\psi$.
By accounting for \eqref{thetapibdd} and the definition \eqref{defArho} of~$\Arho$ once more, we obtain
\Bsist
  && \frac 12 \iO |\psi(t)|^2
  + \nu \intQt |\nabla\psi|^2
  + \intQt (\xi-\xistar) (-\psi)
  + \rho \intQt \bigl( \sigma-(-\zeta) \bigr) (-\psi)
  \non
  \\
  && \quad {}
  + \intQt \bigl( \pi(\phistar+\chi) - \pi(\phistar-w) \bigr) (-\psi)
  \non
  \\
  && = \intQt \bigl( - \gamma\theta - \nu\Delta\phistar + \xistar - \Arho + \pi(\phistar-w) \bigr) \psi
  \non
  \\
  && \leq \intQt \bigl( \gamma\Mrho + \nu \norma{\Delta\phistar}_\infty + \norma\xistar_\infty - \Arho + \Mpi \bigr) \psi
  = 0 .
  \non
\Esist
Also in this case, the integrals on the \lhs\ involving $\xi$ and $\sigma$ are nonnegative:
indeed, where $\psi\not=0$, we have $\psi>0$ and $\chi<-w$,
whence $\phistar+\chi<\phistar-w\leq\phistar$, so that $\xi\leq\xistar$ and $\sigma\leq-\zeta$.
On the other hand, we~have
\Beq
  \intQt \bigl( \pi(\phistar+\chi) - \pi(\phistar-w) \bigr) (-\psi)
  \geq - L \intQt |\chi+w| \psi
  = - L \intQt |\psi|^2 .
  \non
\Eeq
Hence, we deduce \eqref{lastgronwall} with the new meaning of~$\psi$ and apply the Gronwall lemma.
We obtain $\psi=0$, i.e., $-\chi\leq w$.
Therefore, we have proved that $|\chi|\leq w$, 
\gianni{and this implies that $\chi(t)=0$ for every $t\in[\Tstar,T]$}.\QED
}%

\Brem
\label{DecreasingC}
As announced in \gianni{Remark~\ref{Simpldyn}},
we can show that the function $\normaH{\chi(\cpto)}$ is strictly decreasing while positive
provided that $\rho$ is large enough,
at least under a reinforcement of assumption~\eqref{smallness}.
Indeed, with the notations of \accorpa{chibddC}{thetabddC}, 
we have to require that
\Beq
  \rho > (\gamma+L) \CS\CO|\Omega|^{7/6} \rho + \gianni\Ctilde
  \label{reinf}
\Eeq
where we have set
\Beq
  \gianni\Ctilde :=
  \gianniold{ \gamma C_7 + L C_6}
  + \nu \norma{\Delta\phistar}_\infty
  + \norma\xistar_\infty
  + L \norma\phistar_\infty
  + |\pi(0)| .
  \non
\Eeq
Notice that \eqref{reinf} is true provided that
\Beq
  (\gamma+L)\CS\CO|\Omega|^{7/6}<1
  \aand
  \rho > \frac {\gianni\Ctilde} {1 - (\gamma+L)\CS\CO|\Omega|^{7/6}} \,.
  \non
\Eeq
We multiply \eqref{secondachi} written at the time $t$ by $\chi(t)$ and integrate over~$\Omega$.
By ignoring some nonnegative terms on the \lhs, 
observing that $\sigma\chi=|\chi|$ by the definition of~$\sign$,
and owing to~\accorpa{chibddC}{thetabddC},
we easily obtain
\Bsist
  && \frac 12 \, \frac d {dt} \, \normaH{\chi(t)}^2
  + \rho \iO |\chi(t)|
  \leq \iO
  \bigl(
    \gamma\theta(t) + \nu\Delta\phistar + \xistar - \pi(\chi(t)+\phistar)
  \bigr) \chi(t)
  \non
  \\
  &&\leq
  \bigl(
    \gamma \norma\theta_\infty
    + \nu \norma{\Delta\phistar}_\infty
    + \norma\xistar_\infty
    + L \norma\chi_\infty
    + L \norma\phistar_\infty
    + |\pi(0)|
  \bigr) \iO |\chi(t)| .
  \non
  \\
  &&\leq
  \bigl(
    (\gamma+L) \CS\CO|\Omega|^{7/6} \rho + \gianni\Ctilde
  \bigr) \iO |\chi(t)| .
  \non
\Esist
Therefore, on account of~\eqref{reinf}, we conclude that $(d/dt)\normaH{\chi(t)}^2<0$
while $\norma{\chi(t)}_1>0$, or equivalently $\normaH{\chi(t)}>0$.
\Erem



\section*{Acknowledgements}

This research activity has been performed in the framework of an
Italian-Romanian  {three-year project on ``Nonlinear partial differential equations (PDE) 
with applications in modeling cell growth, chemotaxis and phase transition'' financed by the Italian CNR and the Romanian Academy.} 
Moreover, 
\revdue{a partial support} of the FP7-IDEAS-ERC-StG \#256872
(EntroPhase) and the
\revdue{project Fondazione Cariplo-Regione Lombardia MEGAsTAR ``Matematica d'Eccellenza in biologia
ed ingegneria come accelleratore di una nuona strateGia per l'ATtRattivit\`a dell'ateneo pavese''}
is gratefully acknowledged by the authors.
The present paper 
also benefits from the support of the UEFISCDI project PNII-ID-PCE-2011-3-0027 for VB and~GM, the MIUR-PRIN Grant 2010A2TFX2 ``Calculus of Variations'' for PC and GG, 
the GNAMPA (Gruppo Nazionale per l'Analisi Matematica, la Probabilit\`a e le loro Applicazioni) \rev{of INdAM} (Istituto Nazionale di Alta Matematica) for PC, GG and~ER. 


\vspace{3truemm}

\Begin{thebibliography}{99}

\bibitem{BFPU08}
G. Bartolini, L. Fridman, A. Pisano, E. Usai (eds.),
``Modern Sliding Mode Control \revdue{Theory New} Perspectives and Applications'',
Lecture Notes in Control and Information Sciences {\bf 375},
Springer, 2008.

\bibitem{Barbu}
V. Barbu,
``Nonlinear \revdue{Differential Equations of Monotone Types in Banach Spaces}'',
Springer, New York, 2010.

\pier{%
\bibitem{Brezis}
H. Brezis,
``Op\'erateurs \revdue{Maximaux Monotones et Semi-groupes de Contractions
dans les Espaces} de Hilbert'',
North-Holland Math. Stud.
{\bf 5},
North-Holland,
Amsterdam,
1973.
}%

\bibitem{BrokSpr} 
{M. Brokate, J. Sprekels,}
``Hysteresis and \revdue{Phase Transitions}'',
Springer, New York, 1996.

\bibitem{Cag}
G.~Caginalp,
An analysis of a phase field model of a free boundary,
\revdue{{\it Arch. Rational Mech. Anal.}}~{\bf 92} (1986), 205-245.

\bibitem{CRS11}
\revdue{M.-B.} Cheng, V. Radisavljevic, W.-C. Su,
Sliding mode boundary control of a parabolic PDE system with parameter variations and boundary uncertainties, 
{\it Automatica J. IFAC} {\bf 47} (2011), 381--387.

\pier{%
\bibitem{CGM}
P. Colli, G. Gilardi, G. Marinoschi\revdue{,}
A~boundary control problem for a possibly singular phase field system with dynamic boundary conditions,
\rev{{\it J. Math. Anal. Appl.} {\bf 434} (2016), 432--463.}
}%

\bibitem{CoGiMaRo}
P. Colli, G. Gilardi, G. Marinoschi, E. Rocca,
Optimal control for a phase field system 
with a possibly singular potential,
\rev{{\it Math. Control Relat. Fields}
{\bf 6} (2016), 95--112.} 

\pier{%
\bibitem{duvaut}
G. Duvaut,
R\' esolution d'un probl\`eme de Stefan (fusion d'un bloc de glace \`a z\'ero 
degr\'e), {\it C. R. Acad. Sci. Paris S�r. A-B} {\bf 276} (1973), A1461--A1463.
}%

\bibitem{EFF06}
C. Edwards, E. Fossas Colet, L. Fridman (eds.),
``Advances in Variable Structure and Sliding Mode Control'',
Lecture Notes in Control and Information Sciences {\bf 334},
Springer-Verlag, 2006.

\bibitem{ES99}C. Edwards, S. Spurgeon,
``Sliding Mode Control: Theory and Applications'',
Taylor and Francis, London, 1999.

\bibitem{EllZheng}
{C.M. Elliott, S. Zheng,} 
Global existence and stability of solutions to the phase-field equations, 
in ``Free \revdue{Boundary Problems}'', Internat. Ser. Numer. Math. {\bf 95}, 46--58, Birkh\"auser
Verlag, Basel, (1990).

\pier{%
\bibitem{fremond}
M. Fr\'emond,
``Non-smooth Thermomechanics'',
Springer-Verlag, Berlin, 2002.}

\bibitem{FMI11}
L. Fridman, J. Moreno, R. Iriarte (eds.),
``Sliding Modes After the First Decade of the 21st Century: State of the Art'',
Lecture Notes in Control and Information Sciences {\bf 412},
Springer, 2011.

\rev{\bibitem{GT}
D. Gilbarg, N.S. Trudinger,
``Elliptic \revdue{Partial Differential Equations of Second Order}'', Springer, New York, 1998.}

\bibitem{GraPetSch}
M. Grasselli, H. Petzeltov\'a, G. Schimperna,
Long time behavior of solutions to the Caginalp system with singular potential, 
{\it Z. Anal. Anwend.} {\bf 25} (2006), {51--72}.

\bibitem{HoffJiang}
K.-H. Hoffmann, L.S. Jiang, 
Optimal control of a phase field model for solidification, {\it Numer. Funct. Anal. Optim.} {{\bf 13} (1992), 11--27}. 

\bibitem{HKKY}
K.-H. Hoffmann, N. Kenmochi, M. Kubo, N. Yamazaki, 
Optimal control problems for models of phase-field type with hysteresis of play operator,
{\it Adv. Math. Sci. Appl.} {\bf 17} (2007), {305--336}.

\bibitem{I76}
U. Itkis,
``Control \revdue{Systems of Variable Structure}'',
Wiley, 1976.

\bibitem{KenmNiez}
N. Kenmochi, M. Niezg\'odka,
Evolution systems of nonlinear variational inequalities arising \rev{in} phase change problems, 
{\it Nonlinear Anal.} {\bf 22} (1994), 1163--1180.

\bibitem{LSU}
O.A. Lady\v zenskaja, V.A. Solonnikov, N.N. Ural'ceva:
``Linear and \revdue{Quasilinear Equations of Parabolic Type}'',
Trans. Amer. Math. \revdue{Soc.~{\bf 23}},
Amer. Math. Soc., Providence, RI,
1968.

\bibitem{Lau}
Ph. Lauren\c cot,
Long-time behaviour for a model of phase-field type, 
{\it Proc. Roy. Soc. Edinburgh Sect.~A} {\bf 126} (1996), 167--185.

\bibitem{Levaggi13}
L. Levaggi, 
Existence of sliding motions for nonlinear evolution equations in Banach spaces, 
{\it \pier{Discrete Contin. Dynam. Systems}},
Supplement 2013, 477--487.

\bibitem{LO02}
\revdue{L. Levaggi,
Infinite dimensional systems' sliding motions,
{\it \pier{Eur. J. Control}}~{\bf 8} (2002), 508--516.}

\bibitem{O83}
Y.V. Orlov,
Application of Lyapunov method in distributed systems,
{\it \pier{Autom.} Remote Control} {\bf 44} (1983), 426--430.

\bibitem{O00}
Y.V. Orlov, 
Discontinuous unit feedback control of uncertain infinitedimensional
systems, 
{\it IEEE \pier{Trans.} Automatic Control} {\bf 45} (2000), 834--843.

\bibitem{OU83}
Y.V. Orlov, V.I.  Utkin,
Use of sliding modes in distributed system control
\revdue{problems,}
{\it \pier{Autom.} Remote Control} {\bf 43} (1983), 1127--1135.

\bibitem{OU87}
Y.V. Orlov, V.I.  Utkin,
Sliding mode control in indefinite-dimensional systems, 
{\it Automatica
J. IFAC} {\bf 23} (1987), 753--757.

\bibitem{OU98}
Y.V. Orlov, V.I.  Utkin, 
{Unit sliding mode control in infinite-dimensional systems. Adaptive
learning and control using sliding modes},  
{\it Appl. Math. Comput. Sci.} {\bf 8} (1998), 7--20.

\bibitem{RS06}
R. Rossi, G. Savar\'e,
Gradient flows of non convex functionals in Hilbert spaces and applications, 
{\it ESAIM Control Optim. Calc. Var.} {\bf 12} (2006), \revdue{564--614.}

\bibitem{Show}
R.E. Showalter,
``Monotone \revdue{Operators} in Banach Space and Nonlinear Partial Differential Equations'', 
\revdue{Math. Surveys Monogr. {\bf 49}, AMS, Providence, RI,} 1991.

\bibitem{Simon}
J. Simon,
{Compact sets in the space $L^p(0,T; B)$},
{\it Ann. Mat. Pura Appl.~(4)\/}
{\bf 146} (1987), 65--96.

\bibitem{Utkin92}
V. Utkin,
``Sliding Modes in Control and Optimization'',
Communications and Control Engineering Series,
Springer-Verlag, 1992.

\bibitem{UGS09}
V. Utkin, J. Guldner, J. Shi,
``Sliding Mode Control in Electro-Mechanical Systems'', 2nd Edition
CRC Press, Automation and Control Engineering, 2009.

\bibitem{XLGK13}
H. Xing, D. Li, C. Gao, Y. Kao,
Delay-independent sliding mode control for a class of quasi-linear parabolic distributed parameter systems with time-varying delay,
{\it J. Franklin Inst.} {\bf 350} (2013), 397--418.

\bibitem{YO99}
K.D. Young, \"U. \"Ozg\"uner (eds.),
``Variable Structure Systems, Sliding Mode and Nonlinear \revdue{Control''},
Springer-Verlag, 1999.

\End{thebibliography}

\End{document}

\bye